%% file: main.tex
\documentclass[12pt,letterpaper]{article}
\usepackage{amssymb, amsthm, amsmath}
\usepackage{tikz}
\usepackage{float}
\usepackage{hyperref}
\usepackage{multicol}
\usepackage[margin=1in]{geometry}
\usepackage{longtable}
\usepackage{soul}

\usetikzlibrary{arrows.meta}
\usetikzlibrary{shapes.multipart}

\input{graphs}
\input{tables}

\DeclareMathOperator{\spec}{spec}
\DeclareMathOperator{\diag}{diag}

\newtheorem{theorem}{Theorem}
\newtheorem{corollary}[theorem]{Corollary}
\newtheorem{lemma}[theorem]{Lemma}
\newtheorem{prop}[theorem]{Proposition}

\newtheorem{observation}{Observation}
\theoremstyle{definition}
\newtheorem{definition}{Definition}

\allowdisplaybreaks

\begin{document}

\title{Ordered multiplicity inverse eigenvalue problem for graphs on six vertices}
\author{John Ahn\footnote{Bowdoin College, Brunswick, ME 04011, USA {\tt jtahn@bowdoin.edu}} \and
Christine Alar\footnote{San Francisco State University, San Francisco, CA 94132, USA {\tt dalar@mail.sfsu.edu}} \and
Beth Bjorkman\footnote{Iowa State University, Ames, IA 50011, USA {\tt\{bjorkman,butler,jmsdg7\}@iastate.edu}} \and
Steve Butler\footnotemark[3] \and
Joshua Carlson\footnotemark[3] \and
Audrey Goodnight\footnote{Agnes Scott College, Decatur, GA 30030, USA {\tt audrey.goodnight@gmail.com}} \and
Haley Knox\footnote{Eastern Connecticut State University, Willimantic, CT 06286, USA {\tt knoxh@my.easternct.edu}} \and
Casandra Monroe\footnote{Princeton University, Princeton, NJ 08544, USA {\tt cdm5@princeton.edu}} \and
Michael C. Wigal\footnote{West Virginia University, Morgantown, WV 26505, USA {\tt mcwigal@mix.wvu.edu}}}

\date{\empty}

\maketitle

\begin{abstract}
For a graph $G$, we associate a family of real symmetric matrices, $\mathcal{S}(G)$, where for any $M \in \mathcal{S}(G)$, the location of the nonzero off-diagonal entries of $M$ are governed by the adjacency structure of $G$. The ordered multiplicity \emph{Inverse Eigenvalue Problem of a Graph (IEPG)} is concerned with finding all attainable ordered lists of eigenvalue multiplicities for matrices in $\mathcal{S}(G)$.

For connected graphs of order six, we offer significant progress on the IEPG, as well as a complete solution to the ordered multiplicity IEPG.  We also show that while $K_{m,n}$ with $\min(m,n)\ge 3$ attains a particular ordered multiplicity list, it cannot do so with arbitrary spectrum.
\end{abstract}

\section{Introduction}\label{sec:intro}

A graph $G$ consists of a vertex set $V(G)$ and an edge set $E(G)$.  Given $G$ with vertices $v_1,\dots, v_n$, a real symmetric matrix $M$ is in $\mathcal{S}(G)$ if for all $i \neq j$, $M_{i,j}=0$ if and only if $v_iv_j\notin E(G)$; there are no restrictions on the diagonal entries.

The \emph{spectrum} of a matrix $M$ is the set of eigenvalues of $M$. Let $\lambda_1 < \cdots < \lambda_k$ be the distinct eigenvalues of $M$ in increasing order, and let $\gamma_i$ be the multiplicity of $\lambda_i$ as an eigenvalue of $M$. Then the \emph{ordered multiplicity list} of $M$ is $(\gamma_1, \ldots, \gamma_k)$.  (With this convention, the spectrum of $M$ is $\{\lambda_1^{(\gamma_1)},\lambda_2^{(\gamma_2)},\ldots,\lambda_k^{(\gamma_k)}\}$.) 

The \emph{Inverse Eigenvalue Problem of a Graph (IEPG)} is stated as follows: given $G$ and a set of numbers  $L =\{\ell_1,\ldots, \ell_n\}$, does there exist a matrix $M\in\mathcal{S}(G)$ with spectrum $L$?  

This problem has been completely resolved through graphs on five vertices (see \cite{BIRS}).  More information about the IEPG can be found in the survey of Hogben \cite{H}.   Because of the difficulty of the IEPG, many relaxations have been considered; previous works have examined inverse inertia (see \cite{inertia}), minimum rank and maximum nullity (see \cite{FH}), and the minimum number of distinct eigenvalues (see \cite{q}).

We consider the \emph{ordered multiplicity inverse eigenvalue problem for graphs}, a slight relaxation of the IEPG: given a graph $G$ and an ordered list of integers $\Gamma = (\gamma_1, \ldots, \gamma_k)$, does there exist a matrix $M\in\mathcal{S}(G)$ that attains $\Gamma$ as its ordered multiplicity list?

\begin{figure}[htp!]
\centering
\input{diagram.tex}
\caption{There are 26 equivalence classes; those in blue have one graph while those in yellow have their full membership given in the boxes.  To determine attainable ordered multiplicity lists for a graph, find its equivalence class in the diagram and take \emph{any} path to $\emptyset$; the multiplicity lists (and reversals) that occur on the edges of the path are the only ones attainable.  The graph $G$ attains all multiplicity lists that $H$ attains if and only if there is a directed path from the class containing $G$ to the class containing $H$; the difference in what is attainable are the multiplicity lists which occur on \emph{any} directed path between them.}
\label{fig:soln}
\end{figure}
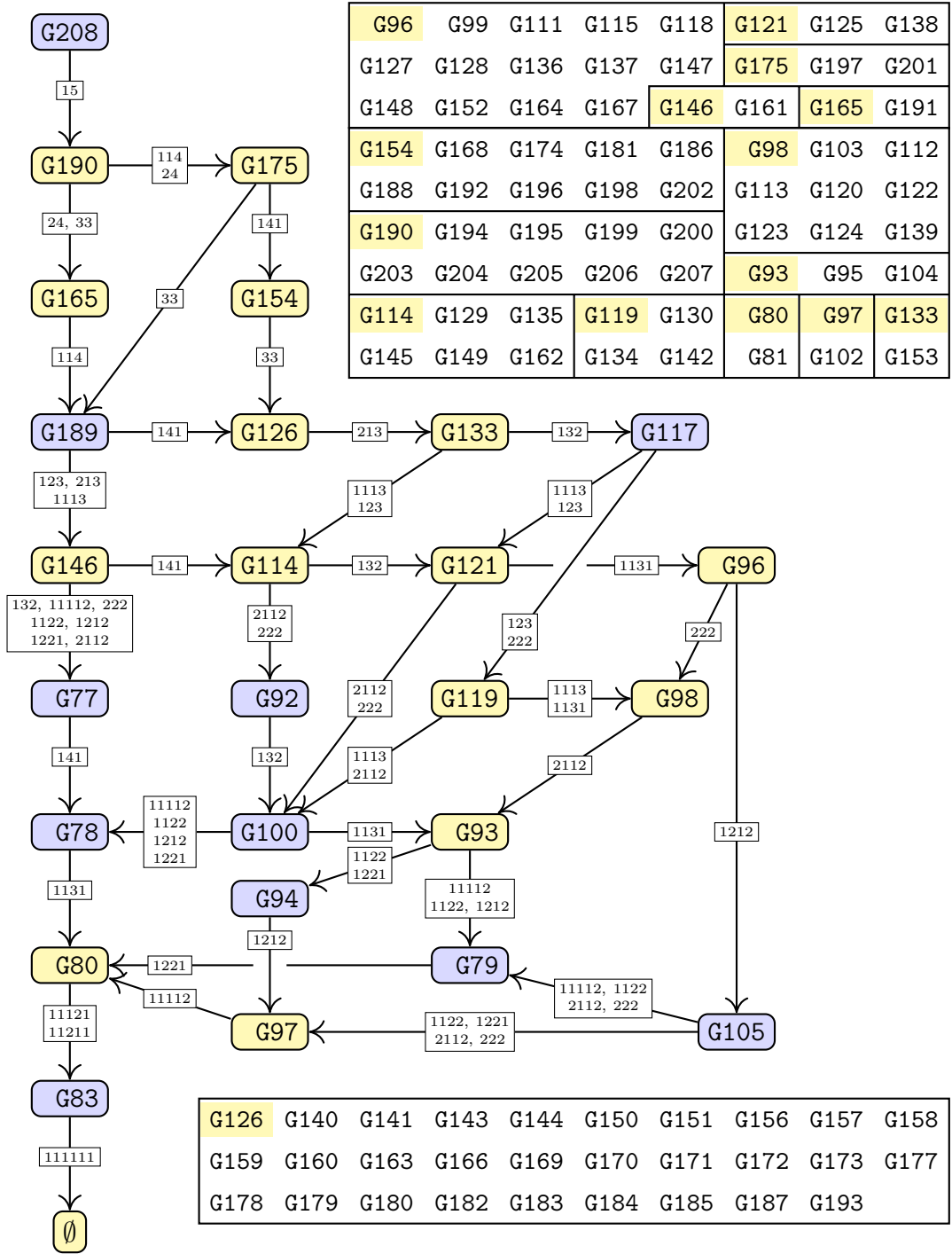

We label graphs using the \emph{Atlas of Graphs} \cite{atlas}; these graphs are also reproduced in the appendix for reference.  In this paper, we solve the ordered multiplicity IEPG for all connected graphs of order six. The result is summarized in Figure~\ref{fig:soln}: the graphs are in 26 different equivalence classes based on what ordered multiplicity lists are attainable.  To determine what lists are attainable, locate the equivalence class it belongs in and then read off all ordered multiplicity lists (and reversals) on the edges of a directed path from that equivalence class to $\emptyset$.  This diagram also gives all possible relationships between equivalence classes, namely the equivalence class containing $G$ attains all ordered multiplicity lists as the equivalence class containing $H$ if and only if there is a directed path from the equivalence class containing $G$ to the one containing $H$.  

We say a graph $G$ is \emph{spectrally arbitrary} for an ordered multiplicity list $(\gamma_1,\ldots,\gamma_k)$ if for any $\lambda_1<\cdots<\lambda_k$, there is a matrix $M\in\mathcal{S}(G)$ with spectrum $\{\lambda_1^{(\gamma_1)},\ldots,\lambda_k^{(\gamma_k)}\}$.  Many of the techniques we use show that a graph is spectrally arbitrary for an ordered multiplicity list. In the appendix, for each ordered multiplicity list, we give all graphs which can attain that list, indicating those which are not known to be attainable with arbitrary spectrum.

We proceed as follows.  In Section~\ref{sec:case5} we will review what is known for the IEPG for graphs on five or fewer vertices which we will build on for the case of six vertices.  In Section~\ref{sec:cloning} we will introduce a technique we call cloning and how it connects with ordered multiplicity lists of eigenvalues.  In Section~\ref{sec:unattain} we justify what ordered multiplicity lists are unattainable, while in Section~\ref{sec:attain} we justify what ordered multiplicity lists are attainable for graphs on six vertices.  In Section~\ref{sec:K33} we show that $K_{m,n}$ with $\min(m,n)\ge3$ is a graph for which the ordered multiplicity IEPG differs from the IEPG.  Finally, in Section~\ref{sec:conclusion}, we give concluding remarks.

Because we are working on graphs with six or fewer vertices, it will be unambiguous to write the ordered multiplicity list $(\gamma_1,\ldots,\gamma_k)$ as $\gamma_1\ldots\gamma_k$. We will say a graph $G$ attains $\gamma_1\ldots\gamma_k$ if there is some $M\in\mathcal{S}(G)$ with multiplicity list $\gamma_1\ldots\gamma_k$; similarly, $G$ does not attain $\gamma_1\ldots\gamma_k$ if there is no $M\in\mathcal{S}(G)$ with multiplicity list $\gamma_1\ldots\gamma_k$.  We note that a graph attains $\gamma_1\ldots\gamma_k$ if and only if it attains $\gamma_k\ldots\gamma_1$, which follows by noting if $M\in\mathcal{S}(G)$ then so also is $-M$.

\section{IEPG for graphs on five or fewer vertices}\label{sec:case5}

The IEPG for all graphs of order at most five was recently solved by Barrett et al.\ \cite{BIRS}.  They showed that for graphs with five or fewer vertices, the IEPG is equivalent to the ordered multiplicity IEPG. Thus, a graph on five or fewer vertices attains a given spectrum if and only if the corresponding multiplicity list is attainable.  Two of the main tools used to solve this problem were the Strong Spectral Property (SSP) and the Strong Multiplicity Property (SMP), introduced in an earlier paper (see \cite{SSP}).

\begin{definition}
An $n\times n$ symmetric matrix $A$ has the \emph{Strong Spectral Property (SSP)} if the only symmetric matrix $X$ satisfying $A \circ X = I \circ X = AX-XA = O$ is $X=O$  (where ``$\circ$'' indicates the Hadamard, or entry-wise, product of matrices).
\end{definition}

\begin{definition}
An $n \times n$ symmetric matrix $A$ satisfies the \emph{Strong Multiplicity Property (SMP)} if the only symmetric matrix $X$ satisfying $A \circ X = I \circ X = AX-XA=O$ and $\text{tr}(A^{i}X) = 0$ for $i = 2, \ldots, n-1$ is $X=O$.
\end{definition}

These properties are important for the IEPG because SSP (SMP) allows us to determine the attainability of certain spectra (ordered multiplicity lists) for many graphs simultaneously.  It should be noted that testing if a matrix has SSP (SMP) reduces to showing if a large linear system has full rank.  This has been implemented and is available online (see \cite{testSSP}); any matrix which is claimed to have SSP follows either from Barrett et al.\ \cite{BIRS} or by using the online implementation.

\begin{theorem}[{Barrett et al.\ \cite{SSP} Theorems 2.10 and 2.20}]\label{thm:SSP1}
If $A\in\mathcal{S}(G)$ has SSP (SMP), then every \emph{supergraph} of $G$ with the same vertex set has an SSP (SMP) realization with the same spectrum (ordered multiplicity list).
\end{theorem}

\begin{theorem}[{Barrett et al.\ \cite{SSP}, Theorem 3.8}]\label{thm:SSP2}
If $A\in\mathcal{S}(G)$ and $B\in\mathcal{S}(H)$ both have SSP (SMP) and $\spec(A)\cap\spec(B)=\emptyset$, then $A\oplus B\in\mathcal{S}(G\cup H)$ has SSP (SMP).
\end{theorem}

Using the above theorems and several constructions, Barrett et al.\ \cite{BIRS} determined that the multiplicity lists given in Table~\ref{tab:5SSP} are attainable with SSP, and those in Table~\ref{tab:5noSSP} are attainable but without SMP or SSP.  Moreover, they established that the graphs can attain any spectrum compatible with the ordered multiplicity lists it attains. Thus, the IEPG and the ordered multiplicity IEPG are equivalent for graphs on five or fewer vertices.

\begin{table}[!th]
\centering
\TABLEFIVESSP
\caption{Realizable ordered multiplicity lists for connected graphs with five or fewer vertices and are attainable with SSP.}
\label{tab:5SSP}
\end{table}

\begin{table}[!th]
\centering
\TABLEFIVENOSSP
\caption{Realizable ordered multiplicity lists for connected graphs with five or fewer vertices which cannot have SMP or SSP.}
\label{tab:5noSSP}
\end{table}

\section{Cloning vertices and ordered multiplicity lists}\label{sec:cloning}

In this section, we introduce \emph{cloning}, a graph operation that, given $G$ and $v \in V(G)$, constructs a new graph $H$ by adding a new vertex $v'$ which is a twin of $v$. This operation is sometimes referred to as duplicating or blow-ups. 

\begin{definition}
Two vertices $u$ and $w$ are \emph{twins} in $H$ if $N_H(u) \setminus \{w\} = N_H(w) \setminus \{u\}$, where $N_H(v)$ is the set of neighbors of $v$ (i.e., vertices which share an edge with $v$).
\end{definition}

Twins do not need to be adjacent. This leads to two variants of cloning: \emph{cloning $v$ with an edge} requires that $v \sim v'$ (i.e., $v$ and $v'$ are adjacent), while \emph{cloning without an edge} requires $v \nsim v'$ (i.e., $v$ and $v'$ are not adjacent).

\begin{theorem}\label{thm:cloning}
Let $G$ be a graph with $M \in \mathcal{S}(G)$ having multiplicity list $(\gamma_1, \ldots , \gamma_i , \ldots, \gamma_k)$ where the eigenvalue $0$ has multiplicity $\gamma_i$.  Then the following two cases are possible:
\begin{enumerate}
\item If the diagonal entry of $M$ corresponding to $v_j$ is zero, then the graph $H$ attained from $G$ by cloning $v_j$ \emph{without} an edge has a matrix $N \in \mathcal{S}(H)$ that attains the multiplicity list $(\gamma_1, \ldots , \gamma_i + 1, \ldots \gamma_k)$.

\item If the diagonal entry of $M$ corresponding to $v_j$ is nonzero, then the graph $H$ attained from $G$ by cloning $v_j$ \emph{with} an edge has a matrix $N \in \mathcal{S}(H)$ that attains the multiplicity list $(\gamma_1, \ldots , \gamma_i + 1, \ldots \gamma_k)$.
\end{enumerate}
\end{theorem}

\begin{proof}
Let $-\lambda_1 \leq \cdots \leq  \lambda_{n - \gamma_i}$ be the nonzero eigenvalues of $M$ where $a$ is the last index such that $-\lambda_a < 0$. Let $\mathbf{x}_1, \ldots, \mathbf{x}_{n - \gamma_i}$ be the corresponding orthonormal eigenvectors. Then 

\begin{align*}
M &= \sum_{k = 1}^{n - \gamma_i} \lambda_k  \mathbf{x}_{k}\mathbf{x}_{k}^T \\
  &=
\begin{pmatrix}\mathbf{x}_1 & \cdots & \mathbf{x}_{n - \gamma_i}\end{pmatrix}
\begin{pmatrix}
-\lambda_1 & \cdots& 0\\
 \vdots& \ddots &\vdots \\
0&\cdots &  \lambda_{n - \gamma_i}
\end{pmatrix}
\begin{pmatrix}
\mathbf{x}_1^T\\
\vdots \\
\mathbf{x}_{n - \gamma_i}^T
\end{pmatrix}\\
&=
\begin{pmatrix} \sqrt{\lambda_1}\mathbf{x}_1 & \cdots & \sqrt{\lambda_{n - \gamma_i}} \mathbf{x}_{n - \gamma_i} \end{pmatrix}
\underbrace{\begin{pmatrix}
- I_a  & O\\
O & I_{n - \gamma_i - a}
\end{pmatrix}}_{=S}
\begin{pmatrix}
\sqrt{\lambda_1}\mathbf{x}_1^T\\
\vdots \\
\sqrt{\lambda_{n - \gamma_i}}\mathbf{x}_{n - \gamma_i}^T
\end{pmatrix}.
\end{align*}
The columns of 
\[
Y = \begin{pmatrix}
\sqrt{\lambda_1}\mathbf{x}_1^T\\
\vdots \\
\sqrt{\lambda_{n - \gamma_i}}\mathbf{x}_{n - \gamma_i}^T
\end{pmatrix}
\]
are an orthogonal representation of $G$ with respect to the indefinite inner product $S$.  That is, if $\mathbf{y}_k$ denotes the $k$-th column then $\mathbf{y}_k^T S \mathbf{y}_\ell = 0$ if and only if $v_k \nsim v_\ell$ in $G$. 
 
Let $Z$ be a $(n - \gamma_i) \times (n + 1)$ matrix with columns as follows,
\begin{itemize}
    \item $\mathbf{z}_k = \mathbf{y}_k$ for $1 \leq k  < j$; 
    \item $\mathbf{z}_j = \mathbf{z}_{j + 1} = \frac{1}{\sqrt{2}}\mathbf{y}_j$;
    \item $\mathbf{z}_k = \mathbf{y}_{k - 1}$ for $j + 1 < k \leq n + 1$. 
\end{itemize}

Now consider the matrix $N=Z^TSZ$. Since $N$ is real symmetric, there exists some graph $H$ such that $N \in\mathcal{S}(H)$. The following two observations will now conclude the proof.



First, use the columns of $Z$ as an orthogonal representation for $H$ with respect to $S$. This corresponds to the graph $G$ with the vertex $v_j$ cloned (i.e.,  columns still have the same orthogonality relationships as given by $Y$).  This will have cloned with an edge if and only if $\mathbf{y}_j^TS\mathbf{y}_j\ne 0$.  The latter holds if and only if the diagonal entry of $M$ corresponding to $v_j$ is nonzero.

Second, the inner product of any two rows of $Z$ agree with the inner product of the corresponding rows of $Y$.  So the nonzero eigenvalues of $N$ are
\[
SZZ^T=\begin{pmatrix}
- I_a  & O\\
O & I_{n - \gamma_i - a}
\end{pmatrix}
\begin{pmatrix}
\lambda_1  & \cdots&0\\
\vdots&\ddots&\vdots\\
0&\cdots& \lambda_{n - \gamma_i}
\end{pmatrix}=
\begin{pmatrix}
-\lambda_1  & \cdots&0\\
\vdots&\ddots&\vdots\\
0&\cdots& \lambda_{n - \gamma_i}
\end{pmatrix},
\]
which are the same as those of $M$.  Hence $N$ has the same spectrum of $M$ with the addition of a single eigenvalue of $0$, giving us the desired ordered multiplicity list.
\end{proof}

\begin{corollary}\label{cor:cloning}
Let $G$ be a graph without isolated vertices, and let $M \in \mathcal{S}(G)$ with multiplicity list $(m_1, \ldots , m_k)$. If a graph $H$ is attained by cloning $v \in V(G)$  with an edge, then $H$ attains the multiplicity list $(m_1+ 1, \ldots, m_k)$.
\end{corollary}

\begin{proof}
By translation we can assume $M$ is positive semi-definite with nullity $m_1$.  If any diagonal entry were $0$, this would force a row and column of zeroes. However, this implies that $G$ contains an isolated vertex, a contradiction. Thus, the entries of the diagonal are nonzero, and so we apply the previous theorem by cloning with an edge.  .
\end{proof}

\section{Unattainable multiplicity lists for graphs}\label{sec:unattain}
In this section we will determine which ordered multiplicity lists are unattainable for connected graphs on six vertices. 

\subsection{Using known graph parameters}
Since we can assign any particular eigenvalue to $0$ by translation, we have the following observations.

\begin{observation}
If $M(G)$ denotes the maximum nullity of a matrix in $\mathcal{S}(G)$, then all entries of a multiplicity list of a matrix in $\mathcal{S}(G)$ are bounded above by $M(G)$.
\end{observation}

\begin{observation}
If $M_{+}(G)$ denotes the maximum nullity of a positive semidefinite matrix in $\mathcal{S}(G)$, then the first (and by reversal from negation, the last) entry of a multiplicity list of a matrix in $\mathcal{S}(G)$ are bounded above by $M_{+}(G)$.
\end{observation}

In general, the computation of $M(G)$ and $M_{+}(G)$ is an open problem. However, for graphs on seven or fewer vertices, known techniques can find find these these values (in particular, the inertia tables---see \cite{atom}). It suffices to provide an upper bound for these parameters, which can be done through the combinatorial parameters $Z(G)$ and $Z_{+}(G)$, respectively known as the zero-forcing number and semidefinite zero-forcing number of a graph. Because of their combinatorial nature, $Z(G)$ and $Z_{+}(G)$ can be easily computed for small graphs through exhaustive analysis. The definition of these parameters, as well as related extensions and results, can be found in the survey of Fallat and Hogben \cite{FH}.  For our purposes, we will use the following result.

\begin{lemma}
[{AIM \cite{z}, Prop. 2.4}]
For any graph $G$, we have $M(G)\le Z(G)$.
\end{lemma}

\begin{lemma}
[Barioli et al.\ \cite{z+}, Theorem 3.5]
For any graph $G$, we have $M_{+}(G)\le Z_{+}(G)$.
\end{lemma}

Another useful parameter is $q(G)$, the minimum number of distinct eigenvalues. 

\begin{observation}
The length of any ordered multiplicity list for $M\in \mathcal{S}(G)$ is at least $q(G)$.
\end{observation}

This is a harder parameter to compute; for connected graphs of order at most six, $q(G)$ was recently determined (see \cite{q}). For all connected graphs of order six, the parameters $Z(G)$, $Z_{+}(G)$, and $q(G)$ are given in Table~\ref{tab:Zdata} and rule out many ordered multiplicity lists.

\begin{table}[ht!]
\TABLEZQ
\caption{The values for $Z(G)$, $Z_+(G)$, and $q(G)$ for connected graphs on six vertices.}
\label{tab:Zdata}
\end{table}

\subsection{Previous results to rule out ordered multiplicity lists}

The following two results, both from \cite{BIRS}, rule out several cases.

\begin{lemma}[{Barrett et al.\ \cite{BIRS}, Lemma 3.3}]\label{lem:unicyclic}
If $G$ is a connected unicyclic graph with odd girth, then at least one of the first or last eigenvalues has multiplicity one.
\end{lemma}

This rules out 2112 and 222 for \texttt{G92}, \texttt{G93}, \texttt{G94}, \texttt{G95}, \texttt{G100}, and \texttt{G104}.

\begin{lemma}[{Barrett et al.\ \cite{BIRS}, Lemmas 2.3 and 5.2}]\label{lem:sun}
A generalized star or a generalized 3-sun does not allow an ordered multiplicity list with consecutive multiple eigenvalues.
\end{lemma}

This rules out 1122, 1221, and 2211 for \texttt{G77} and \texttt{G78} (generalized stars) and \texttt{G94} (3-sun).

The inverse eigenvalue problem for cycles has been determined by Fernandes and Fonseca \cite{FF}, and in particular it follows that 1212 and 2121 are not attainable for the graph $C_6$ (\texttt{G105}).

\subsection{Remaining Cases}

\begin{figure}
\centering
\begin{tabular}{c@{\qquad}c@{\qquad}c@{\qquad}c@{\qquad}c@{\qquad}c}

\begin{tikzpicture}[auto, vertex/.style={circle, draw, inner sep=2pt, thick},scale = 0.35]       
    \node (v1) at (2.5,-2) [vertex] {};
    \node (v2) at (0,3) [vertex] {};
    \node (v3) at (0,-1) [vertex] {};
    \node (v4) at (2.5,1) [vertex] {};
    \node (v5) at (5,-1) [vertex] {};
    \node (v6) at (5,3) [vertex] {};
    \draw [thick] (v1) -- (v4);
    \draw [thick] (v2) -- (v3);
    \draw [thick] (v2) -- (v4);
    \draw [thick] (v3) -- (v4);
    \draw [thick] (v4) -- (v5);
    \draw [thick] (v4) -- (v6);  
    \draw [thick] (v5) -- (v6);
    
    \end{tikzpicture}
 
    &

\begin{tikzpicture}[auto, vertex/.style={circle, draw, inner sep=2pt, thick}, rotate=90,scale = 0.35]       
    \node (v1) at (0,2) [vertex] {};
    \node (v2) at (2,4) [vertex] {};
    \node (v3) at (2,2) [vertex] {};
    \node (v4) at (2,0) [vertex] {};
    \node (v5) at (4,2) [vertex] {};
    \node (v6) at (6,2) [vertex] {};
    \draw [thick] (v1) -- (v2);
    \draw [thick] (v1) -- (v3);
    \draw [thick] (v3) -- (v5);
    \draw [thick] (v1) -- (v4);
    \draw [thick] (v4) -- (v5);
    \draw [thick] (v5) -- (v2);
    \draw [thick] (v5) -- (v6);
    \end{tikzpicture}
 &
   
\begin{tikzpicture}[auto, vertex/.style={circle, draw, inner sep=2pt, thick}, scale = 0.35]       
    \node (v1) at (2,6) [vertex] {};
    \node (v2) at (2,4) [vertex] {};
    \node (v3) at (1,2) [vertex] {};
    \node (v4) at (3,2) [vertex] {};
    \node (v5) at (0,0) [vertex] {};
    \node (v6) at (4,0) [vertex] {};
    \draw [thick] (v1) -- (v2);
    \draw [thick] (v2) -- (v3);
    \draw [thick] (v2) -- (v4);
    \draw [thick] (v3) -- (v5);
    \draw [thick] (v4) -- (v6);
    \draw [thick] (v4) -- (v5);
    \draw [thick] (v3) -- (v6);
    \end{tikzpicture}
    &
    
    \begin{tikzpicture}[auto, vertex/.style={circle, draw, inner sep=2pt, thick}, scale = 0.35,rotate=90]       
    \node (v1) at (0,2) [vertex] {};
    \node (v2) at (2,4) [vertex] {};
    \node (v3) at (2,0) [vertex] {};
    \node (v4) at (4,2) [vertex] {};
    \node (v5) at (5,0) [vertex] {};
    \node (v6) at (5,4) [vertex] {};
    \draw [thick] (v1) -- (v2);
    \draw [thick] (v1) -- (v3);
    \draw [thick] (v1) -- (v4);
    \draw [thick] (v3) -- (v2);
    \draw [thick] (v4) -- (v5);
    \draw [thick] (v4) -- (v6);  
    \draw [thick] (v4) -- (v2);
    \draw [thick] (v4) -- (v3);
    \end{tikzpicture}
    
    &
    
    \begin{tikzpicture}[auto, vertex/.style={circle, draw, inner sep=2pt, thick}, scale = 0.35]       
    \node (v1) at (2,6) [vertex] {};
    \node (v2) at (2,4) [vertex] {};
    \node (v3) at (1,2) [vertex] {};
    \node (v4) at (3,2) [vertex] {};
    \node (v5) at (0,0) [vertex] {};
    \node (v6) at (4,0) [vertex] {};
    \draw [thick] (v1) -- (v2);
    \draw [thick] (v2) -- (v3);
    \draw [thick] (v2) -- (v4);
    \draw [thick] (v3) -- (v5);
    \draw [thick] (v4) -- (v6);
    \draw [thick] (v4) -- (v5);
    \draw [thick] (v3) -- (v6);
    \draw [thick] (v3) -- (v4);
    \end{tikzpicture}
  &

     \begin{tikzpicture}[auto, vertex/.style={circle, draw, inner sep=2pt, thick}, scale = 0.35]       
    \node (v1) at (0,4) [vertex] {};
    \node (v2) at (2,5) [vertex] {};
    \node (v3) at (4,4) [vertex] {};
    \node (v4) at (2,3) [vertex] {};
    \node (v5) at (0,1) [vertex] {};
    \node (v6) at (4,1) [vertex] {};
    \draw [thick] (v1) -- (v2);
    \draw [thick] (v1) -- (v4);
    \draw [thick] (v1) -- (v5);
    \draw [thick] (v3) -- (v2);
    \draw [thick] (v3) -- (v4);
    \draw [thick] (v3) -- (v6);
    \draw [thick] (v5) -- (v6);
    \draw [thick] (v2) -- (v4);
    \end{tikzpicture}
    \\
    \texttt{G117} &
    \texttt{G121} &
    \texttt{G125} &
    \texttt{G133} &
    \texttt{G138} &
    \texttt{G153} 
    \\
    132, 231 & 132, 231 & 132, 231 & 312, 213 & 132, 231 & 312, 213\\
    312, 213
    
\end{tabular}
\caption{The remaining unattainable ordered multiplicity lists.}
\label{fig:remaining}
\end{figure}
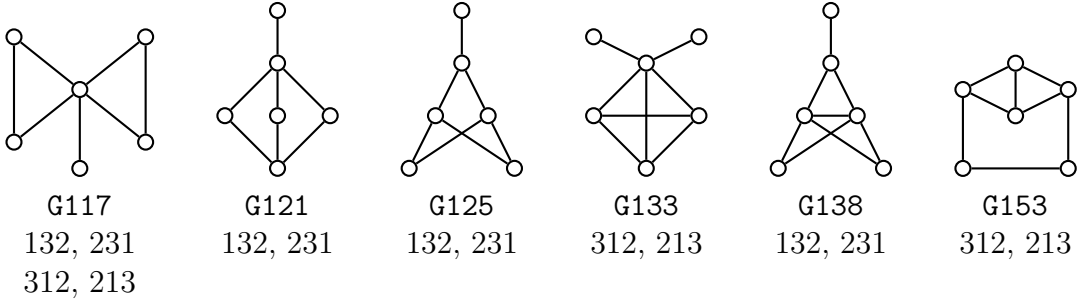

After applying the graph parameters, Lemma \ref{lem:unicyclic}, and Lemma \ref{lem:sun}, the remaining unattainable cases are shown in Figure~\ref{fig:remaining}.

Proving that these cases are unattainable will be done by contradiction; namely, by assuming the existence of a matrix that achieves the specified ordered multiplicity list on a graph. An examination of the corresponding orthogonal representation allows us to argue that two of the vectors are scalar multiples. This pair of vectors allows us to ``declone'' the original graph to a smaller graph with a corresponding matrix that has an impossible ordered multiplicity list.

Recall that cloning works by taking a single vertex and creating a pair of vertices where the vector corresponding with each vertex is a scalar multiple of the original.  Decloning does the opposite: if we can show that the vectors corresponding to a pair of vertices must be scalar multiples, then the two vertices must be twins. We can then delete one of the twins and decrease an entry in the ordered multiplicity list.

\begin{lemma}[Decloning]\label{lem:declone}
Let $G$ be a graph with $M=Q^TSQ\in\mathcal{S}(G)$ having multiplicity list $(\gamma_1,\ldots,\gamma_k)$ where the eigenvalue $0$ has multiplicity $\gamma_i$.  Further assume that $S$ is symmetric and has dimension $n-\gamma_i$, where $n$ is the order of $G$. If the columns of $Q$ corresponding to vertices $a$ and $b$ are scalar multiples, then there exists $N\in\mathcal{S}(H)$ where $H$ is attained from $G$ by deleting vertex $a$ and $N$ has ordered multiplicity list $(\gamma_1,\ldots,\gamma_i-1,\ldots,\gamma_k)$.
\end{lemma}

\begin{proof}
Let $Q=(\mathbf{x}_1\,\mathbf{x}_2\,\cdots\,\mathbf{x}_n)$, and without loss of generality assume that $x_1$ and $x_2$ are scalar multiples, i.e., $\mathbf{x}_2=\alpha\mathbf{x}_1$. Let $\widehat{Q}=(\sqrt{1+\alpha^2}\mathbf{x}_2\,\cdots\,\mathbf{x}_n)$; we claim $N=\widehat{Q}^TS\widehat{Q}$.

Note that two vertices $u$ and $v$ in the graph are not adjacent if and only if $\mathbf{x}_u^TS\mathbf{x}_v=0$. Since scaling by a nonzero value does not change this, we have that the adjacencies are preserved and $N\in\mathcal{S}(H)$.  It remains to check the spectrum of $N$, but for this we note that $QQ^T=\widehat{Q}\widehat{Q}^T$ since the dot products of any pair of corresponding rows are identical.  Since the nonzero portion of the spectrum comes from $SQQ^T=S\widehat{Q}\widehat{Q}^T$, the result follows.
\end{proof}

\begin{prop}\label{prop:312}
If a symmetric matrix $M$ has nullity three and ordered multiplicity 312, then $M=Q^TIQ$ where $I$ is the $3\times3$ identity matrix.
\end{prop}
\begin{proof}
The matrix $M$ has spectrum $\{0^{(3)},\lambda^{(1)},\mu^{(2)}\}$. Let $\mathbf{x}$ be a unit eigenvector for $\lambda$ and $\mathbf{y},\mathbf{z}$ be orthonormal eigenvectors for $\mu$. We have that
\[
M=\begin{pmatrix}\mathbf{x}&\mathbf{y}&\mathbf{z}\end{pmatrix}\begin{pmatrix}\lambda&0&0\\0&\mu&0\\0&0&\mu\end{pmatrix}\begin{pmatrix}\mathbf{x}^T\\\mathbf{y}^T\\\mathbf{z}^T\end{pmatrix}
=\begin{pmatrix}\sqrt\lambda\mathbf{x}&\sqrt\mu\mathbf{y}&\sqrt\mu\mathbf{z}\end{pmatrix}\begin{pmatrix}1&0&0\\0&1&0\\0&0&1\end{pmatrix}\underbrace{\begin{pmatrix}\sqrt\lambda\mathbf{x}^T\\\sqrt\mu\mathbf{y}^T\\\sqrt\mu\mathbf{z}^T\end{pmatrix}}_{=Q}.\qedhere
\]
\end{proof}

\begin{theorem}
The graphs \emph{\texttt{G117}} and \emph{\texttt{G133}} cannot attain ordered multiplicity lists \emph{312} or \emph{213}.
\end{theorem}

\begin{proof}
Label the vertices of \texttt{G117} as shown in Figure~\ref{fig:G117}.  Suppose $M\in\mathcal{S}(\texttt{G117})$ attains ordered multiplicity list 312 with nullity three.  Applying Proposition~\ref{prop:312}, we can write $M=Q^TIQ$, where $Q$ is a matrix which forms an orthogonal representation for \texttt{G117}.

\begin{figure}[ht!]
\centering
\begin{tikzpicture}[auto, vertex/.style={circle, draw, inner sep=2pt, thick},scale = 0.4]       
    \node (v1) at (2.5,-2) [vertex,label=left:$v_1$] {};
    \node (v2) at (0,3) [vertex,label=left:$v_4$] {};
    \node (v3) at (0,-1) [vertex,label=left:$v_3$] {};
    \node (v4) at (2.5,1) [vertex,label=above:$v_2$] {};
    \node (v5) at (5,-1) [vertex,label=right:$v_5$] {};
    \node (v6) at (5,3) [vertex,label=right:$v_6$] {};
    \draw [thick] (v1) -- (v4);
    \draw [thick] (v2) -- (v3);
    \draw [thick] (v2) -- (v4);
    \draw [thick] (v3) -- (v4);
    \draw [thick] (v4) -- (v5);
    \draw [thick] (v4) -- (v6);  
    \draw [thick] (v5) -- (v6);
    
\end{tikzpicture}

\caption{The graph \texttt{G117}.}
\label{fig:G117}
\end{figure}

Since $v_1$ and $v_3$ are not adjacent, the corresponding vectors $\mathbf{v}_1$ and $\mathbf{v}_3$ are orthogonal and thus form a plane in $\mathbb{R}^3$.  Both $v_5$ and $v_6$ are not adjacent to $v_1$ and $v_3$, so $\mathbf{v}_5$ and $\mathbf{v}_6$ are orthogonal to this plane and must be scalar multiples.

The decloning lemma implies that the graph with vertex $e$ deleted has a matrix that attains 212.  However, this graph is odd unicyclic and thus cannot attain 212 by Lemma~\ref{lem:unicyclic}.

A similar argument establishes the result for \texttt{G133}.
\end{proof}

\begin{theorem} 
The graph \emph{\texttt{G153}} cannot attain ordered multiplicity lists \emph{213} or \emph{312}.
\end{theorem}

\begin{proof}
Label the vertices of \texttt{G153} as shown in Figure~\ref{fig:G153}.
Suppose that $M \in \mathcal{S}(\texttt{G153})$ attains ordered multiplicity list 312 with nullity three.  Applying Proposition~\ref{prop:312} we can write $M=Q^TIQ$, where $Q$ forms an orthogonal representation for \texttt{G153}.

\begin{figure}[ht!]
\centering
\begin{tikzpicture}[auto, vertex/.style={circle, draw, inner sep=2pt, thick},scale=0.9]       
    \node (v1) at (0,0) [vertex, label=below:$v_1$] {};
    \node (v2) at (0,1)[vertex, label=above:$v_2$] {};
    \node (v3) at (1,0.5) [vertex, label=right:$v_3$] {};
    \node (v4) at (1,-1) [vertex, label=right:$v_4$] {};
    \node (v5) at (-1,-1) [vertex, label=left:$v_5$] {};
    \node (v6) at (-1,0.5) [vertex, label=left:$v_6$] {};
    \draw [thick] (v1) -- (v2);
    \draw [thick] (v1) -- (v3);
    \draw [thick] (v1) -- (v6);
    \draw [thick] (v2) -- (v3);
    \draw [thick] (v2) -- (v6);
    \draw [thick] (v3) -- (v4);
    \draw [thick] (v4) -- (v5);
    \draw [thick] (v5) -- (v6);
\end{tikzpicture}
\caption{The graph \texttt{G153}.}
\label{fig:G153}
\end{figure}

If $\mathbf{v}_1$ and $\mathbf{v}_2$ are scalar multiples, the decloning lemma implies a matrix for $C_5$ with multiplicity list $212$, which is impossible by Lemma~\ref{lem:unicyclic}.

If $\mathbf{v}_1$ and $\mathbf{v}_2$ are not scalar multiples, then $\mathbf{v}_1$ and $\mathbf{v}_2$ form a plane in $\mathbb{R}^3$. Since vertices $v_4$ and $v_5$ are not adjacent to $v_1$ and $v_2$, the vectors $\mathbf{v}_4$ and $\mathbf{v}_5$ are orthogonal to the aforementioned plane. Thus, $\mathbf{v}_4$ and $\mathbf{v}_5$ are scalar multiples, an impossibility given that $v_4$ and $v_5$ have distinct neighbors.

In either case, we get a contradiction; thus, \texttt{G153} cannot attain 312.
\end{proof}

For the remaining cases, we first establish an analogous result to Proposition~\ref{prop:312}.

\begin{prop}\label{prop:132}
If a symmetric matrix $M$ has nullity three and ordered multiplicity \emph{132}, then $M=Q^TSQ$ where $S=\diag(-1,1,1)$.
\end{prop}
\begin{proof}
The matrix $M$ has spectrum $\{-\lambda^{(1)},0^{(3)},\mu^{(2)}\}$. Let $\mathbf{x}$ be a unit eigenvector for $-\lambda$ and $\mathbf{y},\mathbf{z}$ be orthonormal eigenvectors for $\mu$. We have
\[
M=\begin{pmatrix}\mathbf{x}&\mathbf{y}&\mathbf{z}\end{pmatrix}\begin{pmatrix}-\lambda&0&0\\0&\mu&0\\0&0&\mu\end{pmatrix}\begin{pmatrix}\mathbf{x}^T\\\mathbf{y}^T\\\mathbf{z}^T\end{pmatrix}
=\begin{pmatrix}\sqrt\lambda\mathbf{x}&\sqrt\mu\mathbf{y}&\sqrt\mu\mathbf{z}\end{pmatrix}\begin{pmatrix}-1&0&0\\0&1&0\\0&0&1\end{pmatrix}\underbrace{\begin{pmatrix}\sqrt\lambda\mathbf{x}^T\\\sqrt\mu\mathbf{y}^T\\\sqrt\mu\mathbf{z}^T\end{pmatrix}}_{=Q}.\qedhere
\]
\end{proof}

\begin{theorem}
The graph \emph{\texttt{G117}} cannot attain ordered multiplicity lists \emph{132} or \emph{231}.
\end{theorem}
\begin{proof}
Label the vertices of \texttt{G117} as shown in Figure~\ref{fig:G117}.  Suppose that $M\in\mathcal{S}(\texttt{G117})$ attains ordered multiplicity list 132 with nullity three.  Applying Proposition~\ref{prop:132}, $M=Q^TSQ$ where $S=\diag(-1,1,1)$.

Since $v_1$, $v_3$, and $v_5$ have distinct sets of neighbors, no two of $\mathbf{v}_1$, $\mathbf{v}_3$ and $\mathbf{v}_5$ are scalar multiples of each other.

Because $S$ is invertible, $R=\big({\mathbf{v}_1^TS\atop\mathbf{v}_3^TS}\big)$ has rank two and nullity one. Since $v_1,v_3\notin N(v_5)$, we have $R\mathbf{v}_5=\big({0\atop 0}\big)$; similarly, $v_1,v_3\notin N(v_6)$, so $R\mathbf{v}_6=\big({0\atop 0}\big)$. Thus, $\mathbf{v}_5$ and $\mathbf{v}_6$ are scalar multiples. A similar argument allows us to conclude that $\mathbf{v}_3$ and $\mathbf{v}_4$ are scalar multiples.

Now we can apply the decloning lemma \emph{twice} (i.e., once for each scalar multiple pair) to produce a matrix for the graph $K_{1,3}$ which attains 112.  But $Z_{+}(K_{1,3})=1$, so the end terms of any multiplicity list of $K_{1,3}$ must be $1$, a contradiction.
\end{proof}

In the following proposition, we introduce a tool that shows two vectors are scalar multiples, a technique similar to the one used in the previous results.

\begin{prop}\label{prop:132perp}
If $S$ is a $3{\times}3$ symmetric invertible matrix and $\mathbf{x},\mathbf{y}$ are vectors that satisfy $\mathbf{x}^TS\mathbf{x}=\mathbf{y}^TS\mathbf{y}=\mathbf{x}^TS\mathbf{y}=0$, then $\mathbf{x}$ and $\mathbf{y}$ are scalar multiples of each other.
\end{prop}
\begin{proof}
Assume $\mathbf{x}$ and $\mathbf{y}$ are not scalar multiples; then $S\mathbf{x}$ and $S\mathbf{y}$ are also not scalar multiples.  This implies the matrix
\[
R=\begin{pmatrix}
\mathbf{x}^TS\\
\mathbf{y}^TS
\end{pmatrix}
\]
has rank two and nullity one. On the other hand, by our hypothesis we have 
\[
R\mathbf{x}=\begin{pmatrix}
\mathbf{x}^TS\mathbf{x}\\
\mathbf{y}^TS\mathbf{x}
\end{pmatrix}=
\begin{pmatrix}0\\0\end{pmatrix}=
\begin{pmatrix}
\mathbf{x}^TS\mathbf{y}\\
\mathbf{y}^TS\mathbf{y}
\end{pmatrix}=R\mathbf{y}
\]
which shows that $R$ has nullity at least two, a contradiction.
\end{proof}

\begin{theorem}
The graph \emph{\texttt{G121}} cannot attain ordered multiplicity lists \emph{132} or \emph{231}.
\end{theorem}
\begin{proof}
Label the vertices of \texttt{G121} as shown in Figure~\ref{fig:G121}(a).  Suppose that $M\in\mathcal{S}(\texttt{G121})$ has an ordered multiplicity list 132 with nullity three.

\begin{figure}[ht!]
\centering
\begin{minipage}[t]{0.45\textwidth}
\centering
\begin{tikzpicture}[auto, vertex/.style={circle, draw, inner sep=2pt, thick},scale=0.9]       
    \node (v1) at (1,3) [vertex, label=left:$v_1$] {};
    \node (v2) at (1,2) [vertex, label=left:$v_2$] {};
    \node (v3) at (1,1) [vertex, label=left:$v_3$] {};
    \node (v4) at (1,0) [vertex, label=left:$v_4$] {};
    \node (v5) at (0,1) [vertex, label=left:$v_5$] {};
    \node (v6) at (2,1) [vertex, label=left:$v_6$] {};
    \draw [thick] (v1) -- (v2);
    \draw [thick] (v2) -- (v3);
    \draw [thick] (v2) -- (v5);
    \draw [thick] (v2) -- (v6);
    \draw [thick] (v3) -- (v4);
    \draw [thick] (v4) -- (v5);
    \draw [thick] (v4) -- (v6);
\end{tikzpicture}

(a) Labeled graph.
\end{minipage}
\hfil
\begin{minipage}[t]{0.45\textwidth}
\raggedright
\begin{center}
\begin{tikzpicture}[auto, vertex/.style={circle, draw, inner sep=2pt, thick},scale=0.9]       
    \node (v1) at (1,3) [vertex,color=blue, fill=blue, label=left:$v_1$] {};
    \node (v2) at (1,2)[vertex, label=left:$v_2$] {};
    \node (v3) at (1,1) [vertex,fill=blue, color=blue, label=left:$v_3$] {};
    \node (v4) at (1,0) [vertex, label=left:$v_4$] {};
    \node (v5) at (0,1) [vertex, label=left:$v_5$] {};
    \node (v6) at (2,1) [vertex,color=blue, fill=blue, label=left:$v_6$] {};
    \node (v7) at (3,1) [vertex, label=right:$v_6'$] {};
    \draw [thick] (v1) -- (v2);
    \draw [thick] (v2) -- (v3);
    \draw [thick] (v2) -- (v5);
    \draw [thick] (v2) -- (v6);
    \draw [thick] (v3) -- (v4);
    \draw [thick] (v4) -- (v5);
    \draw [thick] (v4) -- (v6);
    \draw [thick] (v4)--(v7)--(v2) (v7)--(v6);
\end{tikzpicture}
\end{center}

(b) A clone of the graph with an edge and marked zero forcing set.
\end{minipage}
\caption{The graph \texttt{G121}.}
\label{fig:G121}
\end{figure}

We claim that $M_{3,3}=M_{5,5}=M_{6,6}=0$.  To see this, suppose that $M_{6,6}\ne 0$; cloning $v_6$ with an edge produces the graph shown in Figure~\ref{fig:G121}(b) which attains an ordered multiplicity list 142.  The set marked in Figure~\ref{fig:G121}(b) is a zero forcing set of order three. Thus, the cloned graph cannot attain 142, a contradiction. Hence, $M_{6,6}=0$, and $M_{3,3}=M_{5,5}=0$ by symmetry.

We apply Proposition~\ref{prop:132} to write $M=Q^TSQ$ where $S = \diag(-1,1,1)$. Since $v_3$, $v_5$, and $v_6$ form an independent set, we have $\mathbf{v}_3^TS\mathbf{v}_6=\mathbf{v}_3^TS\mathbf{v}_5=\mathbf{v}_5^TS\mathbf{v}_6=0$.  Moreover, because $M_{3,3}=M_{5,5}=M_{6,6}=0$, we have $\mathbf{v}_3^TS\mathbf{v}_3=\mathbf{v}_5^TS\mathbf{v}_5=\mathbf{v}_6^TS\mathbf{v}_6=0$.  From Proposition~\ref{prop:132perp}, it follows that $\mathbf{v}_3$, $\mathbf{v}_5$, and $\mathbf{v}_6$ are pairwise scalar multiples.

Applying the decloning lemma \emph{twice} produces a matrix for $P_4$ which attains 112. However $P_4$ can only attain 1111, a contradiction.
\end{proof}

\begin{theorem}
The graphs \emph{\texttt{G125}} and \emph{\texttt{G138}} cannot attain ordered multiplicity lists \emph{132} or \emph{231}.
\end{theorem}

\begin{proof}
Label the vertices of \texttt{G125} as shown in Figure~\ref{fig:G125}(a).  Suppose $M\in\mathcal{S}(\texttt{G125})$ has spectrum $\{-\lambda,0^{(3)},2^{(2)}\}$ with $\lambda>0$ (by scale and shift, this holds without loss of generality).

\begin{figure}[ht!]
\centering
\begin{minipage}[t]{0.45\textwidth}
\centering
 \begin{tikzpicture}[auto, vertex/.style={circle, draw, inner sep=2pt, thick},scale=0.9]       
    \node (v1) at (1.5,3) [vertex, label=right:$v_1$] {};
    \node (v2) at (1.5,2)[vertex, label=right:$v_2$] {};
    \node (v3) at (1,1) [vertex, label=left:$v_3$] {};
    \node (v4) at (2,1) [vertex, label=right:$v_4$] {};
    \node (v5) at (.5,0) [vertex, label=left:$v_5$] {};
    \node (v6) at (2.5,0) [vertex, label=right:$v_6$] {};
    \draw [thick] (v1) -- (v2);
    \draw [thick] (v2) -- (v3);
    \draw [thick] (v2) -- (v4);
    \draw [thick] (v3) -- (v5);
    \draw [thick] (v4) -- (v6);
    \draw [thick] (v4) -- (v5);
    \draw [thick] (v3) -- (v6);
\end{tikzpicture}

(a) Labeled graph.
\end{minipage}\hfil
\begin{minipage}[t]{0.45\textwidth}
\raggedright
\begin{center}
\begin{tikzpicture}[auto, vertex/.style={circle, draw, inner sep=2pt, thick},scale=0.9]       
    \node (v1) at (1.5,3) [vertex, fill=blue, color=blue, label=right:$v_1$] {};
    \node (v2) at (1.5,2)[vertex, label=right:$v_2$] {};
    \node (v3) at (1,1) [vertex, fill=blue, color=blue, label=left:$v_3$] {};
    \node (v4) at (2,1) [vertex, label=right:$v_4$] {};
    \node (v5) at (.5,0) [vertex, fill=blue, color=blue, label=left:$v_5$] {};
    \node (v6) at (2.5,0) [vertex, label=right:$v_6$] {};
    \node (u) at  (1.25, 0) [vertex, label=right:$v_5'$] {};
    \draw [thick] (v1) -- (v2);
    \draw [thick] (v2) -- (v3);
    \draw [thick] (v2) -- (v4);
    \draw [thick] (v3) -- (v5);
    \draw [thick] (v4) -- (v6);
    \draw [thick] (v4) -- (v5);
    \draw [thick] (v3) -- (v6);
    \draw [thick] (v3) -- (u);
    \draw [thick] (v4) -- (u);
    \draw [thick] (v5) -- (u);
\end{tikzpicture}
\end{center}

(b) A clone of the graph with an edge and marked zero forcing set.
\end{minipage}
\caption{The graph \texttt{G125}.}
\label{fig:G125}
\end{figure}

We claim that $M_{5,5}=M_{6,6}=0$.  To see this, suppose that $M_{5,5}\ne0$; cloning $v_5$ with an edge produces the graph shown in Figure~\ref{fig:G125}(b) which attains an ordered multiplicity list 142.  The set marked in Figure~\ref{fig:G125}(b) is a zero forcing set of order three.  Thus, the cloned graph cannot attain 142, a contradiction.  Hence, $M_{5,5}=0$, and $M_{6,6}=0$ by symmetry.

Applying Proposition~\ref{prop:132} gives $M=Q^TSQ$, where $S=\diag(-1,1,1)$. Note $v_5 \nsim v_6$, so $\mathbf{v}_5^TS\mathbf{v}_6=0$. Moreover, $M_{5,5}=M_{6,6}=0$, so $\mathbf{v}_5^TS\mathbf{v}_5=\mathbf{v}_6^TS\mathbf{v}_6=0$. From Proposition~\ref{prop:132perp}, it follows that $\mathbf{v}_5$ and $\mathbf{v}_6$ are scalar multiples.

Applying the decloning lemma, we attain a matrix $N$ for the banner graph as labeled in Figure~\ref{fig:banner}, where $N$ has eigenvalues $\{-\lambda,0^{(2)},2^{(1)}\}$ and $N_{5,5}=0$; particularly,  $N$ has the following form, where $c_i$ are nonzero and $d_i$ are arbitrary: \newline

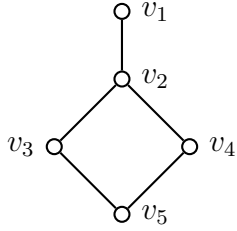
\begin{figure}[ht!]
\centering
\begin{tikzpicture}[auto, vertex/.style={circle, draw, inner sep=2pt, thick},scale=0.9]       
    \node (v1) at (1.5,3) [vertex, label=right:$v_1$] {};
    \node (v2) at (1.5,2)[vertex, label=right:$v_2$] {};
    \node (v3) at (0.5,1) [vertex, label=left:$v_3$] {};
    \node (v4) at (2.5,1) [vertex, label=right:$v_4$] {};
    \node (v5) at (1.5,0) [vertex, label=right:$v_5$] {};
    \draw [thick] (v1) -- (v2);
    \draw [thick] (v2) -- (v3);
    \draw [thick] (v2) -- (v4);
    \draw [thick] (v3) -- (v5);
    \draw [thick] (v4) -- (v5);
\end{tikzpicture}
\caption{The banner graph.}
\label{fig:banner}
\end{figure}

\[
N=\begin{pmatrix}
d_1&c_1&0&0&0\\
c_1&d_2&c_2&c_3&0\\
0&c_2&d_3&0&c_4\\
0&c_3&0&d_4&c_5\\
0&0&c_4&c_5&0
\end{pmatrix}.
\]

The matrix $N$ has rank three; moreover, the first three rows are linearly independent. Therefore the fifth row is a linear combination of the first three rows, which implies that $c_4=\alpha c_2$ and $c_5=\alpha c_3$ for some $\alpha\ne 0$.

Let $R=(N-I)^2$ which has spectrum $\{1^{(4)},(1+\lambda)^2\}$.  The only graphs that attain the ordered multiplicity list 41 are unions of complete graphs with isolated vertices. Since $R_{1,3}=c_1c_2\ne 0$ and $R_{1,5}=0$, we must have that $v_1$ is in the clique and $v_5$ is isolated. In particular
\[
0=R_{2,5}=c_2c_4 + c_3c_5 =\alpha\big(c_2^2+c_3^2)\ne0,
\]
which is impossible.

A similar argument establishes the result for \texttt{G138}.
\end{proof}

\section{Attainable multiplicity lists for graphs}\label{sec:attain}

In this section we establish which ordered multiplicity lists are attainable for connected graphs on six vertices.  We first utilize techniques which will also achieve arbitrary spectrum, and then describe those that fail to do so.

Before we begin, we note that a few special cases have already been done in the literature; we refer the reader elsewhere for details on the following cases.
\begin{itemize}

\item The inverse eigenvalue problem for cycles has been determined by Fernandes and Fonseca.

\begin{theorem}[{Fernandes and Fonseca \cite{FF}, Theorem 3.3}]
The numbers $\lambda_1\le\lambda_2\le\lambda_3\le\cdots\le\lambda_n$ are the spectrum for some matrix $M\in\mathcal{S}(C_n)$ if and only if
\[
\lambda_1\le\lambda_2<\lambda_3\le\lambda_4<\lambda_5\le\lambda_6<\cdots\phantom{.}
\]
or
\[
\lambda_1<\lambda_2\le\lambda_3<\lambda_4\le\lambda_5<\lambda_6<\cdots.
\]
\end{theorem}

(In particular \texttt{G105} attains 11112, 21111, 11121, 12111, 11211, 1122, 2211, 1221, and 2112 spectrally arbitrary.)
\item All trees on six vertices are generalized stars or double stars for which the IEPG has been solved (see \cite{BF}).  (In particular \texttt{G77} attains 1131 and \texttt{G79} attains 1221 with arbitrary spectra.)

\item The graph \texttt{G129} attains 222; a construction can be found in \cite{q}.
\end{itemize}

\subsection{Using SSP for connected graphs of order at most five}

Let $G$ be a disconnected graph on six vertices. From each of its connected components, select an attainable ordered multiplicity list (see Table~\ref{tab:5SSP}). Let $\gamma_1\ldots\gamma_k$ be the ordered multiplicity list built by interlacing in some way these lists. By Theorem~\ref{thm:SSP2}, if each attains its multiplicity list with SSP, $G$ attains $\gamma_1\ldots\gamma_k$ with SSP. By Theorem~\ref{thm:SSP1}, any supergraph of $G$---namely, any connected supergraph $H$ of order six---attains $\gamma_1\ldots\gamma_k$ (with SSP as well). 

Since the components are graphs on at most five vertices, all of the attainable multiplicity lists are spectrally arbitrary. Since SSP preserves spectrum, any ordered multiplicity list constructed in this fashion must be spectrally arbitrary as well. Exhaustively performing the process detailed above gives the results listed in Table~\ref{tab:UseFiveSS}.

\begin{table}[ht!]
\centering
\TABLEA
\caption{Using SSP properties for graphs of order at most five to attain (spectrally arbitrary) ordered multiplicity lists for graphs of order six.}
\label{tab:UseFiveSS}
\end{table}

\subsection{Using cloning for connected graphs of order five}

Let $G$ be a connected graph on five vertices, and let $H$ be the graph constructed by cloning $v \in v(G)$ with an edge. If $G$ attains the ordered multiplicity list $(\gamma_1, \ldots, \gamma_k)$, by Corollary~\ref{cor:cloning} it follows that $H$ attains both $(\gamma_1 + 1, \ldots, \gamma_k)$ and $(\gamma_1, \ldots, \gamma_k + 1)$. 

The method above produces the results listed in Table~\ref{tab:cloning} (note we exclude information that follows from previous results).

\begin{table}[ht!]
\centering
\TABLECLONE
\caption{Using cloning for graphs on five vertices to attain ordered multiplicity lists for graphs on six vertices.}
\label{tab:cloning}
\end{table}

The full application of Theorem~\ref{thm:cloning} generally requires a constructed matrix; then, confirming the value of the appropriate diagonal entry allows for the manipulation of the interior entries of the corresponding ordered multiplicity list. However, an explicit construction is not always necessary. Given a prescribed multiplicity list, for some graphs $G$  we can sometimes guarantee that a particular diagonal entry of any $M \in \mathcal{S}(G)$ must be zero (or nonzero), as we now demonstrate.

\begin{prop}
Let $M\in\mathcal{S}(\emph{\texttt{G40}})$ have nullity two, then the diagonal entry of $M$ corresponding to the leaf is nonzero.  In particular, for any ordered multiplicity list of \emph{\texttt{G40}} with a $2$ we can clone the leaf vertex with an edge to get \emph{\texttt{G144}} and change the $2$ to a $3$.
\end{prop}
\begin{proof}
In Figure~\ref{fig:cloning} we have \texttt{G40} and the two possible graphs that result from cloning without an edge (\texttt{G111}) and with an edge (\texttt{G144}).  In addition we have marked minimal zero forcing sets for the two clones.

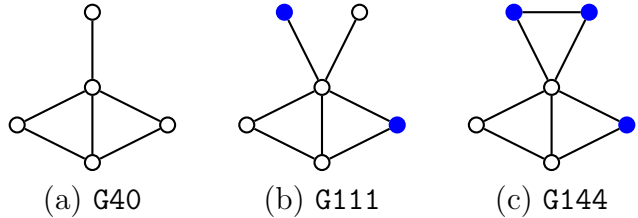
\begin{figure}[ht!]
\centering
\begin{tabular}{c@{\qquad}c@{\qquad}c}
\begin{tikzpicture}
\node[circle, draw, inner sep=2pt, thick] (a) at (0,2) {};
\node[circle, draw, inner sep=2pt, thick] (b) at (0,1) {};
\node[circle, draw, inner sep=2pt, thick] (c) at (-1,0.5) {};
\node[circle, draw, inner sep=2pt, thick] (d) at (0,0) {};
\node[circle, draw, inner sep=2pt, thick] (e) at (1,0.5) {};
\draw[thick] (a)--(b)--(c)--(d)--(e)--(b)--(d);
\end{tikzpicture}
&
\begin{tikzpicture}
\node[circle, draw, inner sep=2pt, thick,fill=blue,color=blue] (a1) at (-0.5,2) {};
\node[circle, draw, inner sep=2pt, thick] (a2) at (0.5,2) {};
\node[circle, draw, inner sep=2pt, thick] (b) at (0,1) {};
\node[circle, draw, inner sep=2pt, thick] (c) at (-1,0.5) {};
\node[circle, draw, inner sep=2pt, thick] (d) at (0,0) {};
\node[circle, draw, inner sep=2pt, thick,fill=blue,color=blue] (e) at (1,0.5) {};
\draw[thick] (a1)--(b) (a2)--(b)--(c)--(d)--(e)--(b)--(d);
\end{tikzpicture}
&
\begin{tikzpicture}
\node[circle, draw, inner sep=2pt, thick,fill=blue,color=blue] (a1) at (-0.5,2) {};
\node[circle, draw, inner sep=2pt, thick,fill=blue,color=blue] (a2) at (0.5,2) {};
\node[circle, draw, inner sep=2pt, thick] (b) at (0,1) {};
\node[circle, draw, inner sep=2pt, thick] (c) at (-1,0.5) {};
\node[circle, draw, inner sep=2pt, thick] (d) at (0,0) {};
\node[circle, draw, inner sep=2pt, thick,fill=blue,color=blue] (e) at (1,0.5) {};
\draw[thick] (b)--(a1)--(a2)--(b)--(c)--(d)--(e)--(b)--(d);
\end{tikzpicture}\\
(a) \texttt{G40}&(b) \texttt{G111}&(c) \texttt{G144}
\end{tabular}
\caption{\texttt{G40} and its two clones via the leaf.  Minimal zero forcing sets have been marked for the clones.}
\label{fig:cloning}
\end{figure}

The graph \texttt{G111} has a zero forcing number of two, which would imply that the maximum nullity (and hence also maximum multiplicity of an eigenvalue) is at most two.  Therefore no matrix associated with \texttt{G111} can have an ordered multiplicity list with an entry which is $\ge3$.  On the other hand that is possible for \texttt{G144}.

We can run cloning with the matrix $M$ which has nullity two and produce a matrix that has nullity three.  This can only be possible if we cloned to \texttt{G144} and not \texttt{G111} which means that the diagonal entry corresponding to the leaf vertex is nonzero.

Finally, since we can always translate any particular eigenvalue to $0$ then for any matrix in $\mathcal{S}(\texttt{G40})$ with an entry of $2$ first we translate so that the entry corresponds to an eigenvalue of $0$, apply the preceding argument, and then translate back.
\end{proof}

A similar argument works for several other cases which are listed in Table~\ref{tab:strongcloning}.

\begin{table}[ht!]
\centering
\TABLECLONEA
\caption{Using cloning for graphs on five vertices to attain ordered multiplicity lists for graphs on six vertices where a middle entry is increased.}
\label{tab:strongcloning}
\end{table}

Because we know that graphs on five vertices attain their ordered multiplicity lists spectrally arbitrarily and cloning preserves the eigenvalues, we can conclude that these lists found by cloning are also spectrally arbitrary.

\subsection{Constructions of graphs which are spectrally arbitrary}

Several cases were handled by finding matrices, usually through the aid of orthogonal representations. Constructions with SSP allowed multiple cases to be handled simultaneously.

\begin{prop}
For \emph{\texttt{G96}} we can attain \emph{222} spectrally arbitrary using SSP matrices.
\end{prop}
\begin{proof}
For $b>0$, the following matrix for \texttt{G96} attains $\{0^{(2)},1^{(2)},(1+4b^2)^{(2)}\}$ with SSP. 
\[
\left(\begin{array}{rrrrrr}
4 b^{2} & \sqrt{2} b & -\sqrt{2} b & 0 & 0 & 0 \\
\sqrt{2} b & 1 & 0 & b & \phantom{-}0 & \phantom{-}0 \\
-\sqrt{2} b & 0 & 1 & b & 0 & 0 \\
0 & b & b & 4 b^{2} & b & b \\
0 & 0 & 0 & b & 1 & 0 \\
0 & 0 & 0 & b & 0 & 1
\end{array}\right)
\]
By choosing $b$, we can make the ratio of the two gaps between the eigenvalues arbitrarily large.  By scaling and shifting, this establishes attainability with arbitrary spectrum.
\end{proof}

By Theorem~\ref{thm:SSP1}, we have that all supergraphs which contain \texttt{G96} can attain the ordered multiplicity list 222 with arbitrary spectrum. Thus, the graphs  \texttt{G111}, \texttt{G114}, \texttt{G118}, \texttt{G121}, \texttt{G135}, \texttt{G136}, \texttt{G137}, \texttt{G140}, \texttt{G145}, \texttt{G146}, \texttt{G147}, \texttt{G148}, \texttt{G149}, \texttt{G157}, \texttt{G158}, \texttt{G159}, \texttt{G161}, \texttt{G162}, \texttt{G163}, \texttt{G164}, \texttt{G166}, \texttt{G167}, \texttt{G171}, \texttt{G173}, \texttt{G180}, \texttt{G182}, \texttt{G184}, \texttt{G186}, \texttt{G187}, \texttt{G188}, \texttt{G196}, \texttt{G197}, \texttt{G198}, and \texttt{G204} are spectrally arbitrary for 222.

\begin{prop}
For \emph{\texttt{G105}} we can attain \emph{222} spectrally arbitrary using SSP matrices.
\end{prop}
\begin{proof}
The matrix
\[
\left(\begin{array}{rrrrrr}
-1 & -1 & 0 & 0 & \phantom{-}0 & a \\
-1 & 0 & a & 0 & 0 & 0 \\
0 & a & a^{2} & a & 0 & 0 \\
0 & 0 & a & -1 & 1 & 0 \\
0 & 0 & 0 & 1 & 0 & a \\
a & 0 & 0 & 0 & a & a^{2}
\end{array}\right)\in\mathcal{S}(\texttt{G105})
\]
has SSP and has spectrum 
\[
\textstyle \big\{(\frac12a^2 - \frac12\sqrt{a^4 + 10a^2 + 5} - \frac12)^{(2)},0^{(2)},(\frac12a^2 + \frac12\sqrt{a^4 + 10a^2 + 5} - \frac12)^{(2)}\big\}.
\]
Since we can scale the spectrum, it suffices to show that the following ratio of the absolute value of the two nonzero eigenvalues contains the interval $[1,\infty)$:
\[
\frac{a^2 + \sqrt{a^4 + 10*a^2 + 5} - 1}{-a^2 + \sqrt{a^4 + 10*a^2 + 5} + 1}.
\]
This is continuous for $a\ge1$ and if $a=1$ we get $1$.  Furthermore, the numerator has growth $a^2$ and the denominator approaches $6$ so the ratio is unbounded.
\end{proof}

By Theorem~\ref{thm:SSP1}, we have that all supergraphs which contain \texttt{G105} can attain the ordered multiplicity list 222 with arbitrary spectrum. Thus, the graphs  \texttt{G127}, \texttt{G128}, \texttt{G151}, \texttt{G152}, \texttt{G154}, \texttt{G174}, and \texttt{G175} are spectrally arbitrary for 222.

\begin{prop}
For \emph{\texttt{G99}} we can attain \emph{222} spectrally arbitrary using non-SSP matrices.
\end{prop}

\begin{proof}
For $a>0$, the following matrix for \texttt{G99} attains $\{0^{(2)},1^{(2)},(1+3a^2)^{(2)}\}$.
\[
\left(\begin{array}{rrrrrr}
1 & a & 0 & \phantom{-}0 & 0 & \phantom{-}0 \\
a & 3 \, a^{2} & a & a & 0 & 0 \\
0 & a & 1 & 0 & -a & 0 \\
0 & a & 0 & 1 & a & 0 \\
0 & 0 & -a & a & 3 \, a^{2} & a \\
0 & 0 & 0 & 0 & a & 1
\end{array}\right)
\]
By choosing $a$ we can make the ratios of the two gaps arbitrarily big. By scaling and shifting, this establishes attainability with arbitrary spectrum. 

It is impossible for any matrix to have 222 and SSP for this graph as \texttt{G112} is a supergraph which cannot attain 222 because $q(\texttt{G112})=4$.
\end{proof}

\begin{prop}
For \emph{\texttt{G189}} we can attain \emph{141} spectrally arbitrary using SSP matrices.
\end{prop}
\begin{proof}
The matrix
\[
\left(\begin{array}{rrrrrr}
0 & \phantom{-}0 & \phantom{-}0 & 1 & 1 & a \\
0 & 0 & 0 & 1 & 1 & a \\
0 & 0 & 0 & 1 & 1 & a \\
1 & 1 & 1 & -1 & -1 & 0 \\
1 & 1 & 1 & -1 & -1 & 0 \\
a & a & a & 0 & 0 & a^{2}
\end{array}\right)\in\mathcal{S}(\texttt{G189})
\]
has SSP and has characteristic polynomial $p(x)=x^4(x^2-(a^2-2)x-(5a^2+6))$.  The spectrum is 
\[
\big\{{\textstyle\frac12}\big((a^2-2)-\sqrt{a^4+16a^2+28}\big),0^{(4)},{\textstyle\frac12}\big((a^2-2)+\sqrt{a^4+16a^2+28}\big)\big\}.
\]
Since we can scale the spectrum, it suffices to show that the following ratio of the absolute value of the two nonzero eigenvalues contains the interval $[1,\infty)$:
\[
\frac{(a^2-2)+\sqrt{a^4+16a^2+28}}{-(a^2-2)+\sqrt{a^4+16a^2+28}}.
\]
This is continuous for $a\ge \sqrt 2$ and if $a=\sqrt{2}$ we get $1$. Furthermore, the numerator has growth $a^2$ and the denominator approaches $10$ so the ratio is unbounded.
\end{proof}

By Theorem~\ref{thm:SSP1}, all supergraphs which contain \texttt{G189} can attain the ordered multiplicity list 141 with arbitrary spectrum.  Thus, the graphs \texttt{G197}, \texttt{G199}, \texttt{G201}, \texttt{G203}, and \texttt{G206} are spectrally arbitrary for 141.

The remaining case for 141 is \texttt{G204}.  Consider
\[
\left(\begin{array}{rrrrrr}
a^{2} & 0 & a & a & \phantom{-}a & \phantom{-}a \\
0 & -1 & -1 & -1 & 1 & 1 \\
a & -1 & 0 & 0 & 2 & 2 \\
a & -1 & 0 & 0 & 2 & 2 \\
a & 1 & 2 & 2 & 0 & 0 \\
a & 1 & 2 & 2 & 0 & 0
\end{array}\right) \in\mathcal{S}(\texttt{G204}),
\]
which has eigenvalues $\{-5, 0^{(4)}, 4+a^2\}$. Arbitrary spectrums are attained through appropriate choices of $a\ge1$, scaling, and shifting.

\begin{prop}
For \emph{\texttt{G151}} we can attain \emph{213}, \emph{312}, \emph{1113}, and \emph{3111} spectrally arbitrary using SSP matrices.
\end{prop}
\begin{proof}
For $a,b>0$, the following matrix for \texttt{G151} has spectrum $\{0^{(3)}, 2a^2, 2a^2+2, a^2+b^2+2\}$ and satisfies SSP:
\[
\left(\begin{array}{rrrrrr}
a^{2} & a & \phantom{-}0 & 0 & a^{2} & 0 \\
a & a^{2} + 2 & a & 0 & 0 & a b \\
0 & a & a^{2} & a^{2} & 0 & 0 \\
0 & 0 & a^{2} & a^{2} + 1 & 1 & -b \\
a^{2} & 0 & 0 & 1 & a^{2} + 1 & -b \\
0 & a b & 0 & -b & -b & b^{2}
\end{array}\right).
\]
Setting $a=b$ gives $312$ with a fixed gap between the last two eigenvalues and an arbitrary gap between the first two; the attainability of an arbitrary spectrum follows. 

For $3111$, by scale and shift we can assume a non-negative spectrum such that $0$ is the eigenvalue of multiplicity $3$, and the gap between the first and second positive eigenvalues is $2$. Set $2a^2$ and $a^2+b^2+2$ to the first and third positive eigenvalues, respectively, and solve for $a$ and $b$. Thus, an arbitrary spectrum is attainable. 
\end{proof}

By Theorem~\ref{thm:SSP1} we have that all supergraphs of \texttt{G151} can attain the ordered multiplicity lists 213, 312, 1113, and 3111 with arbitrary spectrum. Thus, the graphs \texttt{G171} and \texttt{G187} are spectrally arbitrary for 213, 312, 1113, and 3111.

\begin{prop}
For \emph{\texttt{G163}} we can attain \emph{213}, \emph{312}, \emph{1113}, and \emph{3111} spectrally arbitrary using SSP matrices.
\end{prop}

\begin{proof}
For $a,b>0$, the following matrix for \texttt{G163} has spectrum $\{0^{(3)}, 1, 1 + 3a^{2}, 1 + a^2 + 2b^{2}\}$ and satisfies SSP:
\[
\left(\begin{array}{rrrrrr}
a^{2} + b^{2} & b^{2} & \phantom{-}b & 0 & a^{2} & \phantom{-}a \\
b^{2} & a^{2} + b^{2} & b & a & -a^{2} & 0 \\
b &  b & 1 & 0 & 0 & 0 \\
0 & a & 0 & 1 & - a & 0 \\
a^{2} & -a^{2} & 0 & -a & 2 \, a^{2} & a \\
a & 0 & 0 & 0 &  a & 1
\end{array}\right).
\]
By appropriate choices of $a$, $b$, scaling, and translating, this attains 312, 213, 3111, and 1113 spectrally arbitrary.
\end{proof}

\subsection{Two distinct eigenvalues}
For the graphs \texttt{G154}, \texttt{G168}, \texttt{G174}, \texttt{G175}, \texttt{G181}, \texttt{G186}, \texttt{G188}, \texttt{G192}, \texttt{G196}, \texttt{G197}, \texttt{G198}, \texttt{G201}, and \texttt{G202}, $q(G)=2$ and $Z_{+}(G)=3$ (see Table \ref{tab:Zdata}), which implies that all these graphs attain 33. (For more on matrix realizations for these graphs, see \cite{q}.)  For the graph \texttt{G204}, this attains 33 by the matrix below on the left; and 42 (24) by the matrix below on the right.
\[
\left(\begin{array}{rrrrrr}
0 & 0 & 1 & 1 & 1 & 1 \\
0 & 0 & 1 & 1 & -1 & -1 \\
1 & 1 & 0 & 0 & 1 & -1 \\
1 & 1 & 0 & 0 & -1 & 1 \\
1 & -1 & 1 & -1 & 0 & 0 \\
1 & -1 & -1 & 1 & 0 & 0
\end{array}\right)
\qquad\qquad
\left(\begin{array}{rrrrrr}
1 & \phantom{-}0 & \phantom{-}1 & -1 & 2 & -1 \\
0 & 1 & 1 & 1 & 1 & 2 \\
1 & 1 & 2 & 0 & 3 & 1 \\
-1 & 1 & 0 & 2 & -1 & 3 \\
2 & 1 & 3 & -1 & 5 & 0 \\
-1 & 2 & 1 & 3 & 0 & 5
\end{array}\right)
\]

When a matrix for a graph has two distinct eigenvalues, we can modify the matrix to get additional attainable ordered multiplicity lists.  The following will suffice for our purposes (generalizations are possible for graphs with larger order).

\begin{lemma}
Let $G \neq K_6$ be a connected graph on six vertices. If $G$ attains ordered multiplicity list \emph{33}, then with arbitrary spectrum it attains multiplicity lists \emph{1113}, \emph{123}, \emph{213}, \emph{312}, \emph{321}, and \emph{3111}.
Similarly, if $G$ attains \emph{42} or \emph{24}, then with arbitrary spectrum it attains \emph{411} and \emph{114}.
\end{lemma}

\begin{proof}
Since $G$ is not the complete graph, there are two non-adjacent vertices, which we assume to be $v_1$ and $v_2$.

Let $M \in \mathcal{S}(G)$ attain the ordered multiplicity list 33 with spectrum $\{0^{(3)},1^{(3)}\}$.  We can write $M=Q^TQ$, where $Q$ is a $3\times6$ matrix whose rows are \emph{any} orthonormal basis of the eigenspace of $1$ (in particular, $M$ is the projection matrix onto the eigenspace associated with eigenvalue $1$).

We claim we can choose our orthonormal basis so that for some $x,y\ne 0$,
\[
Q = \left(\begin{array}{c}x\\0\\0\end{array}\begin{array}{c}0\\y\\0\end{array}\boxed{\begin{array}{cccc}*&*&*&*\\ *&*&*&*\\ *&*&*&*\end{array}}\,\,\right).
\]
To see this, first consider the vectors $\mathbf{a},\mathbf{b},\mathbf{c}$, which form a basis for our eigenspace.
\begin{itemize}
\item Note that no fixed entry can be $0$ for all three of $\mathbf{a},\mathbf{b},\mathbf{c}$; this would imply that the corresponding vertex is isolated, a contradiction given that our graph is connected. Thus, at least one vector has a nonzero first entry. By taking linear combinations, we can assume that the first entry is $0$ for $\mathbf{a}$ and $\mathbf{b}$ and nonzero for $\mathbf{c}$.
\item Run Gram-Schmidt on $\mathbf{a},\mathbf{b},\mathbf{c}$ (in this order) to get an orthonormal set $\mathbf{a}',\mathbf{b}',\mathbf{c}'$, where the first entries of $\mathbf{a}',\mathbf{b}'$ are $0$ and the first entry of $\mathbf{c}'$ is nonzero.  Note the second entry of $\mathbf{c}'$ must be zero; otherwise, using this set as an orthonormal basis for $Q$ would force $M_{1,2} \neq 0$, a contradiction given that $v_1 \nsim v_2$.
\item Now repeat the argument for $\mathbf{a}',\mathbf{b'}$ by taking a linear combination so that the second entry of $\mathbf{a}'$ is zero. Run Gram-Schmidt again to produce $\mathbf{a}'',\mathbf{b}''$.
\item Thus, the rows of $Q$ are (from top to bottom) $\mathbf{c}', \mathbf{b}'', \mathbf{a}''$.
\end{itemize}

Note we can now introduce parameters $\lambda, \mu>0$ to give
\[
\widehat{Q} = \left(\begin{array}{c}\lambda x\\0\\0\end{array}\begin{array}{c}0\\\mu y\\0\end{array}\boxed{\begin{array}{cccc}*&*&*&*\\ *&*&*&*\\ *&*&*&*\end{array}}\,\,\right).
\]

Since the orthogonality of the columns of $\widehat{Q}$ agree with the orthogonality of the columns of $Q$, the matrix $\widehat{M}=\widehat{Q}^T\widehat{Q}\in\mathcal{S}(G)$. Because the nonzero eigenvalues of $\widehat{M}$ are the norms of the rows, $\widehat{M}$ has spectrum $\{0^{(3)},1,1+(\lambda^2-1)x^2,1+(\mu^2-1)y^2 \}$. Appropriate choices of $\lambda$ and $\mu$, combined with scaling and translation, arbitrarily attain any spectrum that starts or ends with an eigenvalue of multiplicity $3$.

A similar argument handles the 42 case.
\end{proof}

This lemma establishes that 1113, 123, 213, 312, 321, 3111 are all attainable with arbitrary spectrum for graphs \texttt{G154}, \texttt{G174}, and \texttt{G175};
similarly, 411 and 114 are attainable with arbitrary spectrum for \texttt{G204}.

\subsection{Graph minor results}
We now turn to results for graphs which attain certain ordered multiplicity lists, but which are not enough to prove we do so arbitrarily.  We start with the following result which connects SSP and graph minors.

\begin{theorem}[{Barrett et al.\ \cite{BIRS}, Theorem 6.12}]\label{thm:minor}
Let $G$ be attained from $H$ by contraction of a single edge, and let $M\in\mathcal{S}(G)$ have SSP and ordered multiplicity list $(\gamma_1,\ldots,\gamma_k)$.  Then there is $
N\in\mathcal{S}(H)$ with ordered multiplicity list $(\gamma_1,\ldots,\gamma_k,1)$.
\end{theorem}

Applying Theorem~\ref{thm:minor} by looking for minors on graphs of order six gives the results listed in Table~\ref{tab:minors} (we exclude information that follows from previous results).

Note Theorem~\ref{thm:minor} does not guarantee spectrally arbitrary results; the newly appended one on the ordered multiplicity list might need to be large (see \cite{BIRS} for more information).

\begin{table}[ht!]
\centering
\TABLEMINOR
\caption{Using graph minors to attain ordered multiplicity lists for graphs of order six.}
\label{tab:minors}
\end{table}

\subsection{Graphs with SSP/SMP}

Table~\ref{tab:otherSSP} lists matrices that attain the corresponding ordered multiplicity list for that graph. All these matrices have either SSP or SMP, which gives the results listed in Table~\ref{tab:minors}.

\begin{table}[!ht]
\centering
\begin{tabular}{|c|c|c|}\hline
Graph&Multiplicity list(s)&Matrix\\ \hline
\texttt{G105}&2112& \small $\displaystyle \left(\begin{array}{rrrrrr}
0 & \phantom{-}1 & \phantom{-}0 & \phantom{-}0 & \phantom{-}0 & -2 \\
1 & 0 & 2 & 0 & 0 & 0 \\
0 & 2 & 0 & 1 & 0 & 0 \\
0 & 0 & 1 & 0 & 2 & 0 \\
0 & 0 & 0 & 2 & 0 & 1 \\
-2 & 0 & 0 & 0 & 1 & 0
\end{array}\right)
$ \\ \hline
\texttt{G125}&222& \small $\displaystyle\left(\begin{array}{rrrrrr}
-1 & 0 & \phantom{-}1 & \phantom{-}1 & 1 & \phantom{-}0 \\
0 & -1 & 1 & 1 & -1 & 0 \\
1 & 1 & 0 & 0 & 0 & 0 \\
1 & 1 & 0 & 0 & 0 & 0 \\
1 & -1 & 0 & 0 & -1 & 1 \\
0 & 0 & 0 & 0 & 1 & 1
\end{array}\right)
$
\\ \hline
\texttt{G129}& 132, 231 & \small $\displaystyle\left(\begin{array}{rrrrrr}
\phantom{-}2 & 2 & \phantom{-}1 & \phantom{-}1 & \phantom{-}0 & 0 \\
2 & 1 & 0 & 0 & 0 & -1 \\
1 & 0 & 0 & 0 & 1 & 0 \\
1 & 0 & 0 & 0 & 1 & 0 \\
0 & 0 & 1 & 1 & 2 & 2 \\
0 & -1 & 0 & 0 & 2 & 1
\end{array}\right)
$\\ \hline
\end{tabular}
\caption{Matrices with SSP or SMP}
\label{tab:otherSSP}
\end{table}
\begin{itemize}
\item Because \texttt{G127} is a supergraph of \texttt{G105}, \texttt{G127} attains 2112.

\item Because \texttt{G138} is a supergraph of \texttt{G125}, \texttt{G138} attains 222.

\item Because \texttt{G145}, \texttt{G149}, and \texttt{G162} are supergraphs of \texttt{G129}, they attain 132 and 231.
\end{itemize}

\subsection{Remaining cases}
Table~\ref{tab:remaining} gives the remaining attainable cases.  The construction of these matrices included exhaustive searches for matrices with simple entries (i.e., 0, $\pm 1$), as well as the use of orthogonal representations with respect to some (possibly indefinite) inner product.

\begin{table}
\center
\begin{tabular}{|l|l|}\hline

\small$\displaystyle\left(\begin{array}{rrrrrr}
1 & \phantom{-}1 & \phantom{-}1 & \phantom{-}0 & \phantom{-}0 & \phantom{-}0 \\
1 & 1 & 1 & 0 & 0 & 0 \\
1 & 1 & 3 & 2 & 2 & 2 \\
0 & 0 & 2 & 2 & 0 & 0 \\
0 & 0 & 2 & 0 & 2 & 0 \\
0 & 0 & 2 & 0 & 0 & 2
\end{array}\right)
\begin{array}{l}
\in\mathcal{S}(\texttt{G92})\\ \\
\quad1131\\
\quad1311
\end{array}$ 

&

\small$\displaystyle\begin{pmatrix}
1&\phantom{-}1&\phantom{-}1&\phantom{-}1&\phantom{-}0&\phantom{-}0\\
1&1&1&0&1&0\\
1&1&1&0&0&1\\
1&0&0&0&0&0\\
0&1&0&0&0&0\\
0&0&1&0&0&0
\end{pmatrix}
\begin{array}{l}
\in\mathcal{S}(\texttt{G94})\\ \\
\quad2121\\
\quad1212
\end{array}$ 

\\ \hline

\small $\displaystyle\left(\begin{array}{rrrrrr}
0 & 0 & \sqrt{2} & 0 & 0 & 0 \\
0 & 0 & \sqrt{2} & 0 & 0 & 0 \\
\sqrt{2} & \sqrt{2} & -1 & 0 & 1 & 1 \\
0 & 0 & 0 & -1 & 1 & -1 \\
0 & 0 & 1 & 1 & -1 & 1 \\
0 & 0 & 1 & -1 & 1 & -1
\end{array}\right)
\begin{array}{l}
\in\mathcal{S}(\texttt{G114})\\ \\
\quad132\\
\quad231
\end{array}$ 

&

\small
$\displaystyle\left(\begin{array}{rrrrrr}
0 & 1 & 0 & \phantom{-}0 & 0 & \phantom{-}0 \\
1 & -1 & 0 & 0 & 1 & 1 \\
0 & 0 & 0 & 1 & -1 & 1 \\
0 & 0 & 1 & 1 & 0 & 0 \\
0 & 1 & -1 & 0 & 0 & 1 \\
0 & 1 & 1 & 0 & 1 & 0
\end{array}\right)
\begin{array}{l}
\in\mathcal{S}(\texttt{G115})\\ \\
\quad222
\end{array}$ 

\\ \hline

\small$\displaystyle\left(\begin{array}{rrrrrr}
1 & 1 & \phantom{-}1 & \phantom{-}1 & \phantom{-}1 & \phantom{-}1 \\
1 & -1 & 0 & 0 & 0 & 0 \\
1 & 0 & 1 & 1 & 0 & 0 \\
1 & 0 & 1 & 1 & 0 & 0 \\
1 & 0 & 0 & 0 & 1 & 1 \\
1 & 0 & 0 & 0 & 1 & 1
\end{array}\right)
\begin{array}{l}
\in\mathcal{S}(\texttt{G117})\\ \\
\quad1131\\
\quad1311
\end{array}$ 

&

\small$\displaystyle\left(\begin{array}{rrrrrr}
0 & \phantom{-}1 & \phantom{-}1 & \phantom{-}1 & \phantom{-}0 & \phantom{-}0 \\
1 & 1 & 1 & 0 & 0 & 0 \\
1 & 1 & 1 & 0 & 0 & 0 \\
1 & 0 & 0 & 0 & 1 & 1 \\
0 & 0 & 0 & 1 & 1 & 1 \\
0 & 0 & 0 & 1 & 1 & 1
\end{array}\right)
\begin{array}{l}
\in\mathcal{S}(\texttt{G130})\\ \\
\quad1131\\
\quad1311
\end{array}$ 

\\ \hline

\small $\displaystyle\left(\begin{array}{rrrrrr}
0 & 0 & 0 & \sqrt{3} & \sqrt{3} & 0 \\
0 & 0 & 0 & \sqrt{3} & \sqrt{3} & 0 \\
0 & 0 & 0 & \sqrt{3} & \sqrt{3} & 0 \\
\sqrt{3} & \sqrt{3} & \sqrt{3} & 0 & 3 & 0 \\
\sqrt{3} & \sqrt{3} & \sqrt{3} & 3 & 2 & 2 \\
0 & 0 & 0 & 0 & 2 & -1
\end{array}\right)
\begin{array}{l}
\in\mathcal{S}(\texttt{G135})\\ \\
\quad132\\
\quad231
\end{array}$ 

&

\small $\displaystyle\left(\begin{array}{rrrrrr}
\phantom{-}0 & \phantom{-}0 & \phantom{-}0 & 0 & 1 & 1 \\
0 & 0 & 0 & 0 & 1 & 1 \\
0 & 0 & 0 & 0 & 1 & 1 \\
0 & 0 & 0 & -2 & 1 & -1 \\
1 & 1 & 1 & 1 & -1 & 0 \\
1 & 1 & 1 & -1 & 0 & -1
\end{array}\right)
\begin{array}{l}
\in\mathcal{S}(\texttt{G146})\\ \\
\quad132\\
\quad231
\end{array}$ 

\\ \hline

\small$\displaystyle\left(\begin{array}{rrrrrr}
1 & \phantom{-}0 & \phantom{-}0 & \phantom{-}0 & \phantom{-}1 & \phantom{-}1 \\
0 & 1 & 1 & 1 & 1 & 0 \\
0 & 1 & 1 & 1 & 1 & 0 \\
0 & 1 & 1 & 0 & 0 & 0 \\
1 & 1 & 1 & 0 & 1 & 1 \\
1 & 0 & 0 & 0 & 1 & 1
\end{array}\right)
\begin{array}{l}
\in\mathcal{S}(\texttt{G150})\\ \\
\quad1131\\
\quad1311
\end{array}$ 

&

\small $\displaystyle\left(\begin{array}{rrrrrr}
\gamma^2 & \gamma & \phantom{-}\gamma
& 0 & 0 & 0 \\
\gamma & 0 & 2 & -1 & 0 & 0 \\
\gamma & 2 & 0 & 1 & 0 & 0 \\
0 & -1 & 1 & 0 & \gamma & \gamma \\
0 & 0 & 0 & \gamma & \gamma^2 & \gamma^2 \\
0 & 0 & 0 & \gamma & \gamma^2 & \gamma^2
\end{array}\right)
\begin{array}{l}
\gamma^2=\sqrt{10}-2\\[-5pt] \\
\in\mathcal{S}(\texttt{G150})\\ \\
\quad132\\
\quad231
\end{array}$ 

\\ \hline

\small$\displaystyle\left(\begin{array}{rrrrrr}
1 & \phantom{-}3 & 1 & \phantom{-}0 & 1 & \phantom{-}2 \\
3 & 5 & 2 & 1 & 0 & 3 \\
1 & 2 & 7 & 4 &-1 & 0 \\
0 & 1 & 4 & 2 & 0 & 0 \\
1 & 0 &-1 & 0 & -1 & 0 \\
2 & 3 & 0 & 0 & 0 & 2
\end{array}\right)
\begin{array}{l}
\in\mathcal{S}(\texttt{G163})\\ \\
\quad1131\\
\quad1311
\end{array}$ 

&

\small$\displaystyle
\left(\begin{array}{rrrrrr}
0 &  1 &  1 &  0 &  1 & \phantom{-}1 \\
1 &  0 & -1 & -1 & -1 & 0 \\
1 & -1 &  1 &  0 &  0 & 0 \\
0 & -1 &  0 &  1 & -1 & 0 \\
1 & -1 &  0 & -1 &  0 &-1 \\
1 &  0 &  0 &  0 & -1 & 1
\end{array}\right)
\begin{array}{l}
\in\mathcal{S}(\texttt{G163})\\ \\
\quad123\\
\quad321
\end{array}$ 

\\ \hline

\small $\displaystyle\left(\begin{array}{rrrrrr}
5 & \sqrt{5} & 5 & 0 & 0 & 0 \\
\sqrt{5} & 2 & \sqrt{5} & \sqrt{5} & \sqrt{5} & 0 \\
5 & \sqrt{5} & 0 & 0 & -5 & -5 \\
0 & \sqrt{5} & 0 & 5 & 5 & 0 \\
0 & \sqrt{5} & -5 & 5 & 0 & -5 \\
0 & 0 & -5 & 0 & -5 & -5
\end{array}\right)
\begin{array}{l}
\in\mathcal{S}(\texttt{G163})\\ \\
\quad132\\
\quad231
\end{array}$ 

& \\ \hline
\end{tabular}
\caption{Remaining cases.}
\label{tab:remaining}
\end{table}

\section{Ordered Multiplicity IEPG differs from IEPG}\label{sec:K33}
Through five vertices, the ordered multiplicity IEPG and the IEPG are equivalent: a graph attains an ordered multiplicity list if and only if it attains that multiplicity list with arbitrary spectrum. However, for graphs of order at least six, this relationship no longer holds.

\begin{theorem}
The complete bipartite graph $K_{m,n}$ where $\min(m,n) \ge 3$ attains the ordered multiplicity list $(1,m+n-2,1)$, but not spectrally arbitrary.
\end{theorem}

\begin{proof}
The spectrum of the adjacency matrix of $K_{m,n}$ is $\{-\sqrt{mn},0^{(m+n-2)},\sqrt{mn}\}$. Thus, $K_{m,n}$ attains  $(1,m+n-2,1)$. Note the gaps between consecutive eigenvalues are equal. We claim that any $M \in \mathcal{S}(K_{m,n})$ that attains $(1,m+n-2,1)$ preserves this relationship. Thus, $K_{m,n}$ cannot attain $(1,m+n-2,1)$ with arbitrary spectrum.

Let $M \in \mathcal{S}(K_{m,n})$ attain multiplicity list $(1,m + n -2,1)$, where after translation we may assume $0$ is the eigenvalue of multiplicity $m + n - 2$.  Let $-\lambda_1 < \lambda_2$ be the nonzero eigenvalues of $M$ with the corresponding orthogonal eigenvectors $\mathbf{x}_1$ and $\mathbf{x}_2$. Note,
\begin{align*}
M 
= -\lambda_1 \mathbf{x}_1\mathbf{x}_1^T + \lambda_2 \mathbf{x}_2\mathbf{x}_2^T 
&=
\begin{pmatrix}
\mathbf{x}_1 & \mathbf{x}_2\\
\end{pmatrix}
\begin{pmatrix}
-\lambda_1 & 0\\
0 & \lambda_2
\end{pmatrix}
\begin{pmatrix}
\mathbf{x}_1^T\\
\mathbf{x}_2^T
\end{pmatrix}\\
&=
\begin{pmatrix}
\sqrt{\lambda_1} \mathbf{x} _1 & \sqrt{ \lambda_2 }\mathbf{x}_2
\end{pmatrix}
\underbrace{\begin{pmatrix}
-1 & 0 \\
0 & 1\\
\end{pmatrix}}_{=S}
\underbrace{\begin{pmatrix}
\sqrt{\lambda_1} \mathbf{x}_1^T\\
\sqrt{\lambda_2} \mathbf{x}_2^T
\end{pmatrix}}_{=Y}
.
\end{align*}

Let $\mathbf{y}_i$ be the $i^{th}$ column of $Y$; thus, $y_i \in \mathbb{R}^2$. Note the association between $y_i$ and $v_i$ of $K_{m,n}$: two vertices $v_i \nsim v_j$ if and only if $\mathbf{y}_i^T S \mathbf{y}_j = 0$. Moreover, $K_{m,n}$ is connected, so $y_i \neq 0$.

Let $a, b \in K_{m,n}$ be nonadjacent. Since $\min(m,n) \geq 3$, there exists $c$ such that $a,b,$ and $c$ are pairwise nonadjacent. If the corresponding vectors $\mathbf{a}$ and $\mathbf{b}$ are not scalar multiples, then the matrix
\[
\begin{pmatrix}
-a_1 & a_2\\
-b_1 & b_2
\end{pmatrix}
\] 
has rank two. However, the adjacency structure of $K_{m,n}$ requires that
\[
\begin{pmatrix}
-a_1 & a_2\\
-b_1 & b_2
\end{pmatrix}
\begin{pmatrix}
c_1\\
c_2
\end{pmatrix}=
\begin{pmatrix}
0\\
0
\end{pmatrix},
\]
which shows the matrix lacks full rank, a contradiction. Thus, $\mathbf{a}$ and $\mathbf{b}$ are scalar multiples.  By symmetry, the vectors associated with pairwise nonadjacent vertices must be scalar multiples. These vectors must also satisfy $\mathbf{a}^TS\mathbf{a}=0$, and thus must have the form $\big({\phantom{\pm}x\atop \pm x}\big)$. Hence, $Y$ is of the form
\[
\begin{pmatrix}
\alpha_1&\alpha_2&\cdots&\alpha_m&\phantom{-}\beta_1&\phantom{-}\beta_2&\cdots&\phantom{-}\beta_n\\
\alpha_1&\alpha_2&\cdots&\alpha_m&-\beta_1&-\beta_2&\cdots&-\beta_n\\
\end{pmatrix},
\]
for appropriate choice of $\alpha_i$ and $\beta_j$, so $M$ has spectrum $\big\{0^{(m+n-2)},\pm2\sqrt{(\sum\alpha_i^2)(\sum\beta_j^2)}\big\}$. Therefore $K_{m,n}$ cannot attain $(1,m+n-2,1)$ with arbitrary spectrum.
\end{proof}

\begin{corollary}
The graph $K_{3,3}$ (\emph{\texttt{G175}}) attains $141$, but not spectrally arbitrary.
\end{corollary}

\section{Conclusion}\label{sec:conclusion}
We have given a complete solution for the ordered multiplicity IEPG for connected graphs on six vertices. Moreover, many of the techniques used for attainability also allow for an arbitrary spectrum, allowing for significant progress on the IEPG for connected graphs on six vertices. In particular, there are 1326 cases of attainability. Among these, 1285 are known to be done with arbitrary spectrum, one does so without, and 40 remain undetermined (see Table~\ref{tab:IEPG}). Finishing these cases, and thus solving the IEPG for graphs on six vertices, is an open problem.

\begin{table}[!ht]
\centering
\TABLEREMAININGIEPG
\caption{The remaining cases for the IEPG for connected graphs on six vertices.}
\label{tab:IEPG}
\end{table}

We also showed that $K_{m,n}$ with $\min(m,n)\ge 3$ has at least one attainable multiplicity list which cannot be attained spectrally arbitrary.  This shows that the ordered multiplicity IEPG and the IEPG differ for graphs on six or more vertices.

A natural next problem to consider is the ordered multiplicity IEPG for connected graphs on seven or more vertices. While many of the techniques implemented in this work can be utilized in that setting (and indeed solves ``most'' of the cases), the IEPG becomes dramatically harder,  Thus, new tools will likely be need to continue to make progress.

The difficulty lies in part with the number of cases involved, both in terms of the number of graphs and the number of potential ordered multiplicity lists.  In addition, one of the most useful tools we had was SMP and SSP which allowed for simultaneous handling of many cases by establishing a result for a graph and all its supergraphs.  This does have some limitations.

\begin{observation}
If $M\in\mathcal{S}(G)$ attains $\gamma_1\ldots\gamma_k$ and $H$ is a supergraph which does not attain $\gamma_1\ldots\gamma_k$, then $M$ does not have SSP or SMP.
\end{observation}

This can be used to explain why the cases given in Table~\ref{tab:5noSSP} do not have SSP, and for graphs on six vertices can be used to show that there are over forty occurrences where a graph attains an ordered multiplicity list but does so without SSP or SMP.  These cases required either cloning or finding some appropriate orthogonal constructions.  As the number of vertices increases the number of cases needing individual attention (and hence difficulty) will rise as well.

\subsection*{Acknowledgments}
The research was conducted at the 2017 REU program held at Iowa State University, which was supported by NSF DMS 1457443.  Steve Butler was partially supported by a grant from the Simons Foundation (\#427264).

\input{appendix}

\end{document}

%% file: tables.tex
\def\TABLEA{\begin{tabular}{|l|l|l|}\hline
Graph(s) & Supergraphs & Multiplicity lists\\ \hline
\texttt{G54} $\cup$ \texttt{G1} & 
\texttt{G191}, \texttt{G200}, \texttt{G205}, \texttt{G207}, \texttt{G208}&
111111, 11112, 11121, 11211, 12111, 21111, \\ 
&& 1122, 1212, 1221, 2112, 2121, 2211,\\ 
&& 1113, 1131, 1311, 3111, 123, 132,\\
&& 213, 231, 312, 321, 114, 141, 411\\ \hline
\texttt{G48} $\cup$ \texttt{G1}&
\texttt{G140}, \texttt{G141}, \texttt{G143}, \texttt{G156}, \texttt{G157},& 111111, 11112, 11121, 11211, 12111, 21111,\\ &
\texttt{G158}, \texttt{G159}, \texttt{G160}, \texttt{G166}, \texttt{G168},& 1122, 1212, 1221, 2112, 2121, 2211,\\ &
\texttt{G170}, \texttt{G172}, \texttt{G173}, \texttt{G177}, \texttt{G178},& 1113, 1131, 1311, 3111, 123, 132,\\ &
\texttt{G179}, \texttt{G180}, \texttt{G181}, \texttt{G182}, \texttt{G183},& 213, 231, 312, 321 \\ &
\texttt{G184}, \texttt{G185}, \texttt{G186}, \texttt{G188}, \texttt{G189},& \\ &
\texttt{G190}, \texttt{G192}, \texttt{G193}, \texttt{G194}, \texttt{G195},& \\ &
\texttt{G196}, \texttt{G197}, \texttt{G198}, \texttt{G199}, \texttt{G201},& \\ &
\texttt{G202}, \texttt{G203}, \texttt{G204}, \texttt{G206}
 & \\ \hline
\texttt{G18} $\cup$ 2\texttt{G1} & \texttt{G133}, \texttt{G134}, \texttt{G142}, \texttt{G165}, \texttt{G169}& 111111, 11112, 11121, 11211, 12111, 21111, \\
&& 1122, 1212, 1221, 2112, 2121, 2211, \\
&&  1113, 1131, 1311, 3111 \\\hline
\texttt{G44} $\cup$ \texttt{G1} & 
\texttt{G121}, \texttt{G125}, \texttt{G135}, \texttt{G138}, \texttt{G146},& 111111, 11112, 11121, 11211, 12111, 21111, \\ &
\texttt{G149}, \texttt{G154}, \texttt{G161}, \texttt{G162}, \texttt{G175}& 1122, 1212, 1221, 2112, 2121, 2211,\\ && 1131, 1311 \\ \hline

\texttt{G7} $\cup$ \texttt{G7} & \texttt{~G96}, \texttt{~G98}, \texttt{~G99}, \texttt{G103}, \texttt{G111}, &  111111, 11112, 11121, 11211, 12111, 21111,\\ 

\texttt{G42} $\cup$ \texttt{G1} & \texttt{G112}, \texttt{G113}, \texttt{G114}, \texttt{G115}, \texttt{G117}, & 1122, 1212, 1221, 2112, 2121, 2211 \\ 

\texttt{G16} $\cup$ 2\texttt{G1} & \texttt{G118}, \texttt{G119}, \texttt{G120}, \texttt{G122}, \texttt{G123}, & \\ 

& \texttt{G124}, \texttt{G126}, \texttt{G128}, \texttt{G129}, \texttt{G130}, & \\ 
& \texttt{G136}, \texttt{G137}, \texttt{G139}, \texttt{G144}, \texttt{G145}, & \\ 
& \texttt{G147}, \texttt{G148}, \texttt{G150}, \texttt{G151}, \texttt{G152}, & \\ 
& \texttt{G153}, \texttt{G163}, \texttt{G164}, \texttt{G167}, \texttt{G171}, & \\ 
& \texttt{G174}, \texttt{G187} & \\ \hline

\texttt{G34} $\cup$ \texttt{G1}&
\texttt{~G92}, \texttt{~G93}, \texttt{~G95}, \texttt{G104}, \texttt{G127} &  111111, 11112, 11121, 11211, 12111, 21111,\\
\texttt{G38} $\cup$ \texttt{G1} &&  1122, 1212, 1221, 2121, 2211\\ \hline

\texttt{G7} $\cup$ 3\texttt{G1} & \texttt{~G94}, \texttt{~G97}, \texttt{G100}, \texttt{G102} &  111111, 11112, 11121, 11211, 12111, 21111 \\ \hline

\texttt{G13} $\cup$ 2\texttt{G1} & \texttt{~G77}, \texttt{~G78}, \texttt{~G79}, \texttt{~G80}, \texttt{~G81} &  111111, 11121, 11211, 12111\\ \hline

6\texttt{G1} & \texttt{~G83}, \texttt{G105} & 111111 \\ \hline
\end{tabular}}

\def\TABLEZQ{\begin{minipage}{0.2\textwidth}
\[
\begin{array}{r|ccc}
G&Z&Z_{+}&q\\ \hline
\texttt{G77}&4&1&3\\
\texttt{G78}&3&1&4\\
\texttt{G79}&2&1&4\\
\texttt{G80}&2&1&5\\
\texttt{G81}&2&1&5\\
\texttt{G83}&1&1&6\\
\texttt{G92}&3&2&3\\
\texttt{G93}&2&2&4\\
\texttt{G94}&2&2&4\\
\texttt{G95}&2&2&4\\
\texttt{G96}&2&2&3\\
\texttt{G97}&2&2&5\\
\texttt{G98}&2&2&4\\
\texttt{G99}&2&2&3\\
\texttt{G100}&3&2&4\\
\texttt{G102}&2&2&5\\
\texttt{G103}&2&2&4\\
\texttt{G104}&2&2&4\\
\texttt{G105}&2&2&3\\
\texttt{G111}&2&2&3\\
\texttt{G112}&2&2&4\\
\texttt{G113}&2&2&4\\
\texttt{G114}&3&2&3\\
\texttt{G115}&2&2&3\\
\texttt{G117}&3&3&3\\
\texttt{G118}&2&2&3\\
\texttt{G119}&3&3&4\\
\texttt{G120}&2&2&4
\end{array}
\]
\end{minipage}
\hfil
\begin{minipage}{0.2\textwidth}
\[
\begin{array}{r|ccc}
G&Z&Z_{+}&q\\ \hline
\texttt{G121}&3&2&3\\
\texttt{G122}&2&2&4\\
\texttt{G123}&2&2&4\\
\texttt{G124}&2&2&4\\
\texttt{G125}&3&2&3\\
\texttt{G126}&3&3&3\\
\texttt{G127}&2&2&3\\
\texttt{G128}&2&2&3\\
\texttt{G129}&3&2&3\\
\texttt{G130}&3&3&4\\
\texttt{G133}&3&3&3\\
\texttt{G134}&3&3&4\\
\texttt{G135}&3&2&3\\
\texttt{G136}&2&2&3\\
\texttt{G137}&2&2&3\\
\texttt{G138}&3&2&3\\
\texttt{G139}&2&2&4\\
\texttt{G140}&3&3&3\\
\texttt{G141}&3&3&3\\
\texttt{G142}&3&3&4\\
\texttt{G143}&3&3&3\\
\texttt{G144}&3&3&3\\
\texttt{G145}&3&2&3\\
\texttt{G146}&4&2&3\\
\texttt{G147}&2&2&3\\
\texttt{G148}&2&2&3\\
\texttt{G149}&3&2&3\\
\texttt{G150}&3&3&3
\end{array}
\]
\end{minipage}
\hfil
\begin{minipage}{0.2\textwidth}
\[
\begin{array}{r|ccc}
G&Z&Z_{+}&q\\ \hline
\texttt{G151}&3&3&3\\
\texttt{G152}&2&2&3\\
\texttt{G153}&3&3&3\\
\texttt{G154}&3&3&2\\
\texttt{G156}&3&3&3\\
\texttt{G157}&3&3&3\\
\texttt{G158}&3&3&3\\
\texttt{G159}&3&3&3\\
\texttt{G160}&3&3&3\\
\texttt{G161}&4&2&3\\
\texttt{G162}&3&2&3\\
\texttt{G163}&3&3&3\\
\texttt{G164}&2&2&3\\
\texttt{G165}&4&4&3\\
\texttt{G166}&3&3&3\\
\texttt{G167}&2&2&3\\
\texttt{G168}&3&3&2\\
\texttt{G169}&3&3&3\\
\texttt{G170}&3&3&3\\
\texttt{G171}&3&3&3\\
\texttt{G172}&3&3&3\\
\texttt{G173}&3&3&3\\
\texttt{G174}&3&3&2\\
\texttt{G175}&4&3&2\\
\texttt{G177}&3&3&3\\
\texttt{G178}&3&3&3\\
\texttt{G179}&3&3&3\\
\texttt{G180}&3&3&3
\end{array}
\]
\end{minipage}
\hfil
\begin{minipage}{0.2\textwidth}
\[
\begin{array}{r|ccc}
G&Z&Z_{+}&q\\ \hline
\texttt{G181}&3&3&2\\
\texttt{G182}&3&3&3\\
\texttt{G183}&3&3&3\\
\texttt{G184}&3&3&3\\
\texttt{G185}&3&3&3\\
\texttt{G186}&3&3&2\\
\texttt{G187}&3&3&3\\
\texttt{G188}&3&3&2\\
\texttt{G189}&4&3&3\\
\texttt{G190}&4&4&2\\
\texttt{G191}&4&4&3\\
\texttt{G192}&3&3&2\\
\texttt{G193}&3&3&3\\
\texttt{G194}&4&4&2\\
\texttt{G195}&4&4&2\\
\texttt{G196}&3&3&2\\
\texttt{G197}&4&3&2\\
\texttt{G198}&3&3&2\\
\texttt{G199}&4&4&2\\
\texttt{G200}&4&4&2\\
\texttt{G201}&4&3&2\\
\texttt{G202}&3&3&2\\
\texttt{G203}&4&4&2\\
\texttt{G204}&4&4&2\\
\texttt{G205}&4&4&2\\
\texttt{G206}&4&4&2\\
\texttt{G207}&4&4&2\\
\texttt{G208}&5&5&2
\end{array}
\]
\end{minipage}}

\def\TABLEFIVESSP{\begin{tabular}{|l|l|}\hline
Graphs&Attainable Ordered multiplicity lists\\ \hline
\texttt{G1} & 1 \\ \hline
\texttt{G3} & 11 \\ \hline
\texttt{G6} & 111 \\ \hline
\texttt{G7} & 111, 12, 21 \\ \hline
\texttt{G14} & 1111 \\ \hline
\texttt{G13} & 1111, 121 \\ \hline
\texttt{G15} & 1111, 121, 112, 211 \\ \hline
\texttt{G16}, \texttt{G17}&1111, 121, 112, 211, 22 \\ \hline
\texttt{G18} & 1111, 121, 112, 211, 22, 13, 31 \\ \hline
\texttt{G31} & 11111 \\ \hline
\texttt{G29}, \texttt{G30} & 11111, 1121, 1211 \\  \hline
\texttt{G35}, \texttt{G36} & 11111, 1121, 1211, 1112, 2111 \\ \hline
\texttt{G34}, \texttt{G38} & 11111, 1121, 1211, 1112, 2111, 122, 221 \\ \hline
\texttt{G37}, \texttt{G40}, \texttt{G41} & 11111, 1121, 1211, 1112, 2111, 122, 221, 212 \\ 
\texttt{G42}, \texttt{G43}, \texttt{G47} & \\ \hline
\texttt{G44}, \texttt{G46} & 11111, 1121, 1211, 1112, 2111, 122, 221, 212, 131 \\ \hline
\texttt{G45}&11111, 1121, 1211, 1112, 2111, 122, 221, 212, 131, 113, 311 \\ \hline
\texttt{G48}, \texttt{G49} & 11111, 1121, 1211, 1112, 2111, 122, 221, 212, 131, 113, 311, 23, 32 \\ 
\texttt{G50}, \texttt{G51} & \\ \hline
\texttt{G52} & 11111, 1121, 1211, 1112, 2111, 122, 221, 212, 131, 113, 311, 23, 32, 14, 41 \\ \hline 
\end{tabular}}

\def\TABLEFIVENOSSP{\begin{tabular}{|l|l|}\hline
Graphs&Attainable Ordered multiplicity lists\\ \hline
\texttt{G29} & 131\\ \hline
\texttt{G42} & 113, 131, 311\\ \hline
\end{tabular}
}

\def\TABLECLONE{\begin{tabular}{|l|l|l|} \hline
Graph&Cloned graph(s)&Multiplicity lists\\ \hline
\texttt{G29} & \texttt{~G92}, \texttt{G161} & 132, 231\\ \hline
\texttt{G30} & \texttt{G100} & 1122, 1212, 2121, 2211\\ \hline
\texttt{G34} & \texttt{G117}, \texttt{G133}, \texttt{G179} & 1113, 123, 222, 3111, 321\\ \hline
\texttt{G35} & \texttt{G119} & 1113, 3111\\ \hline
\texttt{G36} & \texttt{G130}, \texttt{G150} & 1113, 3111\\ \hline
\texttt{G37} & \texttt{G126}, \texttt{G141}, \texttt{G143}, \texttt{G168} & 1113, 123, 213, 222, 3111, 312, 321\\ \hline
\texttt{G38} & \texttt{G153} & 1113, 123, 222, 3111, 321\\ \hline
\texttt{G40} & \texttt{G144}, \texttt{G156}, \texttt{G177}, \texttt{G192} & 1113, 123, 213, 222, 3111, 312, 321\\ \hline
\texttt{G41} & \texttt{G150}, \texttt{G160}, \texttt{G178}, \texttt{G181} & 1113, 123, 213, 222, 3111, 312, 321\\ \hline
\texttt{G42} & \texttt{G165}, \texttt{G195} & 114, 123, 132, 213, 222, 231, 312, 321, 411\\ \hline
\texttt{G43} & \texttt{G169}, \texttt{G172}, \texttt{G185} & 123, 213, 222, 312, 321\\ \hline
\texttt{G44} & \texttt{G170}, \texttt{G189} & 222\\ \hline
\texttt{G45} & \texttt{G165}, \texttt{G191}, \texttt{G200} & 114, 123, 132, 213, 222, 231, 312, 321, 411\\ \hline
\texttt{G46} & \texttt{G179}, \texttt{G201} & 222\\ \hline
\texttt{G47} & \texttt{G183}, \texttt{G193}, \texttt{G202} & 222\\ \hline
\texttt{G48} & \texttt{G190}, \texttt{G194}, \texttt{G199} & 114, 222, 24, 33, 411, 42\\ \hline
\texttt{G49} & \texttt{G195}, \texttt{G200}, \texttt{G205} & 114, 222, 24, 33, 411, 42\\ \hline
\texttt{G50} & \texttt{G203}, \texttt{G206} & 114, 222, 24, 33, 411, 42\\ \hline
\texttt{G51} & \texttt{G205}, \texttt{G207} & 222, 24, 33, 42\\ \hline
\texttt{G52} & \texttt{G208} & 15, 222, 24, 33, 42, 51\\ \hline
\end{tabular}}

\def\TABLECLONEA{\begin{tabular}{|l|l|l|}\hline
Graph & Cloned graph & Multiplicity lists\\ \hline
\texttt{G29} & \texttt{~G77} & 141 \\ \hline
\texttt{G30} & \texttt{~G78} & 1131, 1311 \\ \hline
\texttt{G30} & \texttt{G100} & 1131, 1311 \\ \hline
\texttt{G30} & \texttt{G114} & 1131, 1311 \\ \hline
\texttt{G34} & \texttt{G133} & 132, 231 \\ \hline
\texttt{G35} & \texttt{G119} & 1131, 1311 \\  \hline
\texttt{G37} & \texttt{G126} & 1131, 1311, 132, 231 \\ \hline
\texttt{G40} & \texttt{G144} & 1131, 1311, 132, 231 \\ \hline
\texttt{G42} & \texttt{G165} & 141 \\ \hline
\texttt{G44} & \texttt{G146} & 141 \\ \hline
\texttt{G46} & \texttt{G161} & 141 \\ \hline
\texttt{G48} & \texttt{G190} & 141 \\ \hline
\texttt{G48} & \texttt{G194} & 141 \\ \hline
\texttt{G49} & \texttt{G195} & 141 \\ \hline
\end{tabular}}

\def\TABLEREMAININGIEPG{\begin{tabular}{|l|l|}\hline
Multiplicity list(s) & Remaining cases\\ \hline
1212, 2121 & \texttt{G94} \\ \hline
1221 & \texttt{G100} \\ \hline
2112 & \texttt{G127}\\ \hline
222 & \texttt{G115}, \texttt{G125}, \texttt{G129}, \texttt{G138} \\ \hline
1131, 1311 & \texttt{G92}, \texttt{G117}, \texttt{G129}, \texttt{G130}, \texttt{G145}, \texttt{G150}, \texttt{G151}, \texttt{G153},  \texttt{G163}, \texttt{G171}, \\ & \texttt{G174}, \texttt{G187}\\ \hline
123, 321 & \texttt{G151}, \texttt{G163}, \texttt{G171}, \texttt{G187}\\ \hline
132, 231 & \texttt{G114}, \texttt{G129}, \texttt{G135}, \texttt{G145}, \texttt{G146}, \texttt{G149}, \texttt{G150}, \texttt{G151}, \texttt{G153}, \texttt{G154}, \\ & \texttt{G162}, \texttt{G163}, \texttt{G169}, \texttt{G171}, \texttt{G174}, \texttt{G175}, \texttt{G187}\\ \hline
\end{tabular}
}

\def\TABLEMINOR{\begin{tabular}{|l|l|l|}\hline
Graph(s)&Minor&Multiplicity lists\\ \hline
\texttt{G100}&\texttt{G34}&1221\\ \hline
\texttt{G129}, \texttt{G145}, \texttt{G151}, \texttt{G153}, &\texttt{G44}&1311, 1131\\
\texttt{G171}, \texttt{G174}, \texttt{G187}&&\\\hline
\texttt{G151}, \texttt{G153}, \texttt{G154}, \texttt{G169}, &\texttt{G48}&321, 123, 231, 132 \\
\texttt{G171}, \texttt{G174}, \texttt{G175}, \texttt{G187}&&\\\hline
\end{tabular}
}

%% file: diagram.tex



\begin{tikzpicture}[every text node part/.style={align=center}]
\node[draw,thick,rounded corners, fill=yellow!35!white] (empty) at (0,0) {$\emptyset$};
\node[draw,thick,rounded corners, fill=blue!15!white]  (G83) at (0,2)  {\texttt{~G83}};
\node[draw,thick,rounded corners, fill=yellow!35!white]  (G80) at (0,4)  {\texttt{~G80}};
\node[draw,thick,rounded corners, fill=blue!15!white]  (G78) at (0,6)  {\texttt{~G78}};
\node[draw,thick,rounded corners, fill=blue!15!white]  (G77) at (0,8)  {\texttt{~G77}};
\node[draw,thick,rounded corners, fill=yellow!35!white] (G146) at (0,10) {\texttt{G146}};
\node[draw,thick,rounded corners, fill=blue!15!white] (G189) at (0,12) {\texttt{G189}};
\node[draw,thick,rounded corners, fill=yellow!35!white] (G165) at (0,14) {\texttt{G165}};
\node[draw,thick,rounded corners, fill=yellow!35!white] (G190) at (0,16) {\texttt{G190}};
\node[draw,thick,rounded corners, fill=blue!15!white] (G208) at (0,18) {\texttt{G208}};
\node[draw,thick,rounded corners, fill=yellow!35!white]  (G97) at (3,3)  {\texttt{~G97}};
\node[draw,thick,rounded corners, fill=blue!15!white] (G100) at (3,6)  {\texttt{G100}};
\node[draw,thick,rounded corners, fill=blue!15!white]  (G92) at (3,8)  {\texttt{~G92}};
\node[draw,thick,rounded corners, fill=yellow!35!white] (G114) at (3,10) {\texttt{G114}};
\node[draw,thick,rounded corners, fill=yellow!35!white] (G126) at (3,12) {\texttt{G126}};
\node[draw,thick,rounded corners, fill=yellow!35!white] (G154) at (3,14) {\texttt{G154}};
\node[draw,thick,rounded corners, fill=yellow!35!white] (G175) at (3,16) {\texttt{G175}};
\node[draw,thick,rounded corners, fill=blue!15!white]  (G94) at (3,5)  {\texttt{~G94}};
\node[draw,thick,rounded corners, fill=blue!15!white]  (G79) at (6,4)  {\texttt{~G79}};
\node[draw,thick,rounded corners, fill=yellow!35!white]  (G93) at (6,6)  {\texttt{~G93}};
\node[draw,thick,rounded corners, fill=yellow!35!white] (G119) at (6,8)  {\texttt{G119}};
\node[draw,thick,rounded corners, fill=yellow!35!white] (G121) at (6,10) {\texttt{G121}};
\node[draw,thick,rounded corners, fill=yellow!35!white] (G133) at (6,12) {\texttt{G133}};
\node[draw,thick,rounded corners, fill=yellow!35!white]  (G98) at (9,8)  {\texttt{~G98}};
\node[draw,thick,rounded corners, fill=yellow!35!white]  (G96) at (10,10) {\texttt{~G96}};
\node[draw,thick,rounded corners, fill=blue!15!white] (G117) at (9,12) {\texttt{G117}};
\node[draw,thick,rounded corners, fill=blue!15!white]  (G105) at (10,3) {\texttt{G105}};

\draw[-{>[scale=2,width=5pt]}, thick] (G208)--(G190);
\draw[-{>[scale=2,width=5pt]}, thick] (G190)--(G165);
\draw[-{>[scale=2,width=5pt]}, thick] (G165)--(G189);
\draw[-{>[scale=2,width=5pt]}, thick] (G189)--(G146);
\draw[-{>[scale=2,width=5pt]}, thick] (G146)--(G77);
\draw[-{>[scale=2,width=5pt]}, thick] (G77)--(G78);
\draw[-{>[scale=2,width=5pt]}, thick] (G78)--(G80);
\draw[-{>[scale=2,width=5pt]}, thick] (G80)--(G83);
\draw[-{>[scale=2,width=5pt]}, thick] (G83)--(empty);
\draw[-{>[scale=2,width=5pt]}, thick] (G190)--(G175);
\draw[-{>[scale=2,width=5pt]}, thick] (G175)--(G189);
\draw[-{>[scale=2,width=5pt]}, thick] (G189)--(G126);
\draw[-{>[scale=2,width=5pt]}, thick] (G146)--(G114);
\draw[-{>[scale=2,width=5pt]}, thick] (G100)--(G78);
\draw[-{>[scale=2,width=5pt]}, thick] (G79)--(G80);
\draw[-{>[scale=2,width=5pt]}, thick] (G97)--(G80);
\draw[-{>[scale=2,width=5pt]}, thick] (G175)--(G154);
\draw[-{>[scale=2,width=5pt]}, thick] (G154)--(G126);
\draw[-{>[scale=2,width=5pt]}, thick] (G114)--(G92);
\draw[-{>[scale=2,width=5pt]}, thick] (G92)--(G100);
\draw[-{>[scale=2,width=5pt]}, thick] (G126)--(G133);
\draw[-{>[scale=2,width=5pt]}, thick] (G133)--(G114);
\draw[-{>[scale=2,width=5pt]}, thick] (G114)--(G121);
\draw[-{>[scale=2,width=5pt]}, thick] (G121)--(G100);
\draw[-{>[scale=2,width=5pt]}, thick] (G100)--(G93);
\draw[-{>[scale=2,width=5pt]}, thick] (G93)--(G79);
\draw[-{>[scale=2,width=5pt]}, thick] (G119)--(G100);
\draw[-{>[scale=2,width=5pt]}, thick] (G93)--(G94);
\draw[-{>[scale=2,width=5pt]}, thick] (G133)--(G117);
\draw[-{>[scale=2,width=5pt]}, thick] (G117)--(G121);
\draw[-{>[scale=2,width=5pt]}, thick] (G121)--(G96);
\draw[-{>[scale=2,width=5pt]}, thick] (G119)--(G98);
\draw[-{>[scale=2,width=5pt]}, thick] (G98)--(G93);
\draw[-{>[scale=2,width=5pt]}, thick] (G96)--(G98);
\draw[-{>[scale=2,width=5pt]}, thick] (G96)--(G105);
\draw[-{>[scale=2,width=5pt]}, thick] (G105)--(G79);
\draw[-{>[scale=2,width=5pt]}, thick] (G105)--(G97);

\draw[color=white, line width=10pt] (7.25,10)--(7.75,10);
\draw[-{>[scale=2,width=5pt]}, thick] (G117)--(G119);

\draw[color=white, line width=10pt] (2.75,4)--(3.25,4);
\draw[-{>[scale=2,width=5pt]}, thick] (G94)--(G97);

\node[draw,fill=white,inner sep=2pt] at (0,1.125) {\tiny 111111};
\node[draw,fill=white,inner sep=2pt] at (0,3.125) {\tiny 11121\\[-7pt] \tiny 11211};
\node[draw,fill=white,inner sep=2pt] at (0,5.125) {\tiny 1131};
\node[draw,fill=white,inner sep=2pt] at (0,7.125) {\tiny 141};
\node[draw,fill=white,inner sep=2pt] at (0,9.125) {\tiny 132, 11112, 222\\[-7pt] \tiny 1122, 1212\\[-7pt] \tiny 1221, 2112};
\node[draw,fill=white,inner sep=2pt] at (0,11.125) {\tiny 123, 213\\[-7pt] \tiny 1113};
\node[draw,fill=white,inner sep=2pt] at (0,13.125) {\tiny 114};
\node[draw,fill=white,inner sep=2pt] at (0,15.125) {\tiny 24, 33};
\node[draw,fill=white,inner sep=2pt] at (0,17.125) {\tiny 15};
\node[draw,fill=white,inner sep=2pt] at (1.5,16) {\tiny 114\\[-7pt] \tiny 24};
\node[draw,fill=white,inner sep=2pt] at (1.5,14) {\tiny 33};
\node[draw,fill=white,inner sep=2pt] at (1.5,12) {\tiny 141};
\node[draw,fill=white,inner sep=2pt] at (3,15.125) {\tiny 141};
\node[draw,fill=white,inner sep=2pt] at (3,13.125) {\tiny 33};
\node[draw,fill=white,inner sep=2pt] at (1.5,10) {\tiny 141};
\node[draw,fill=white,inner sep=2pt] at (3,9.125) {\tiny 2112\\[-7pt] \tiny 222};
\node[draw,fill=white,inner sep=2pt] at (3,7.125) {\tiny 132};
\node[draw,fill=white,inner sep=2pt] at (3,4.375) {\tiny 1212};
\node[draw,fill=white,inner sep=2pt] at (1.5,6) {\tiny 11112\\[-7pt] \tiny 1122\\[-7pt] \tiny 1212\\[-7pt] \tiny 1221};
\node[draw,fill=white,inner sep=2pt] at (1.5,4) {\tiny 1221};
\node[draw,fill=white,inner sep=2pt] at (1.5,3.5) {\tiny 11112};
\node[draw,fill=white,inner sep=2pt] at (4.5,6) {\tiny 1131};
\node[draw,fill=white,inner sep=2pt] at (6,5) {\tiny 11112\\[-7pt] \tiny 1122, 1212};
\node[draw,fill=white,inner sep=2pt] at (4.5,12) {\tiny 213};
\node[draw,fill=white,inner sep=2pt] at (4.5,11) {\tiny 1113\\[-7pt] \tiny 123};
\node[draw,fill=white,inner sep=2pt] at (4.5,10) {\tiny 132};
\node[draw,fill=white,inner sep=2pt] at (4.5,8) {\tiny 2112\\[-7pt] \tiny 222};
\node[draw,fill=white,inner sep=2pt] at (4.5,7) {\tiny 1113\\[-7pt] \tiny 2112};
\node[draw,fill=white,inner sep=2pt] at (7.5,12) {\tiny 132};
\node[draw,fill=white,inner sep=2pt] at (7.5,11) {\tiny 1113\\[-7pt] \tiny 123};
\node[draw,fill=white,inner sep=2pt] at (8.5,10) {\tiny 1131};
\node[draw,fill=white,inner sep=2pt] at (4.5,5.5) {\tiny 1122\\[-7pt] \tiny 1221};
\node[draw,fill=white,inner sep=2pt] at (7.5,8) {\tiny 1113\\[-7pt] \tiny 1131};
\node[draw,fill=white,inner sep=2pt] at (7.5,7) {\tiny 2112};
\node[draw,fill=white,inner sep=2pt] at (9.5,9) {\tiny 222};
\node[draw,fill=white,inner sep=2pt] at (6.75,9) {\tiny 123\\[-7pt] \tiny 222};
\node[draw,fill=white,inner sep=2pt] at (10,6) {\tiny 1212};
\node[draw,fill=white,inner sep=2pt] at (6,3) {\tiny 1122, 1221\\[-7pt] \tiny 2112, 222};
\node[draw,fill=white,inner sep=2pt] at (8,3.5) {\tiny 11112, 1122\\[-7pt] \tiny 2112, 222};

\node[fill=yellow!35!white] at ( 4.750,18.125) {\small\texttt{~G96}};
\node at ( 5.875,18.125) {\small\texttt{~G99}};
\node at ( 7.000,18.125) {\small\texttt{G111}};
\node at ( 8.125,18.125) {\small\texttt{G115}};
\node at ( 9.250,18.125) {\small\texttt{G118}};
\node[fill=yellow!35!white] at (10.375,18.125) {\small\texttt{G121}};
\node at (11.500,18.125) {\small\texttt{G125}};
\node at (12.625,18.125) {\small\texttt{G138}};

\node at ( 4.750,17.500) {\small\texttt{G127}};
\node at ( 5.875,17.500) {\small\texttt{G128}};
\node at ( 7.000,17.500) {\small\texttt{G136}};
\node at ( 8.125,17.500) {\small\texttt{G137}};
\node at ( 9.250,17.500) {\small\texttt{G147}};
\node[fill=yellow!35!white] at (10.375,17.500) {\small\texttt{G175}};
\node at (11.500,17.500) {\small\texttt{G197}};
\node at (12.625,17.500) {\small\texttt{G201}};

\node at ( 4.750,16.875) {\small\texttt{G148}};
\node at ( 5.875,16.875) {\small\texttt{G152}};
\node at ( 7.000,16.875) {\small\texttt{G164}};
\node at ( 8.125,16.875) {\small\texttt{G167}};
\node[fill=yellow!35!white] at ( 9.250,16.875) {\small\texttt{G146}};
\node at (10.375,16.875) {\small\texttt{G161}};
\node[fill=yellow!35!white] at (11.500,16.875) {\small\texttt{G165}};
\node at (12.625,16.875) {\small\texttt{G191}};

\node[fill=yellow!35!white] at ( 4.750,16.250) {\small\texttt{G154}};
\node at ( 5.875,16.250) {\small\texttt{G168}};
\node at ( 7.000,16.250) {\small\texttt{G174}};
\node at ( 8.125,16.250) {\small\texttt{G181}};
\node at ( 9.250,16.250) {\small\texttt{G186}};
\node[fill=yellow!35!white] at (10.375,16.250) {\small\texttt{~G98}};
\node at (11.500,16.250) {\small\texttt{G103}};
\node at (12.625,16.250) {\small\texttt{G112}};
\node at ( 4.750,15.625) {\small\texttt{G188}};
\node at ( 5.875,15.625) {\small\texttt{G192}};
\node at ( 7.000,15.625) {\small\texttt{G196}};
\node at ( 8.125,15.625) {\small\texttt{G198}};
\node at ( 9.250,15.625) {\small\texttt{G202}};
\node at (10.375,15.625) {\small\texttt{G113}};
\node at (11.500,15.625) {\small\texttt{G120}};
\node at (12.625,15.625) {\small\texttt{G122}};
\node[fill=yellow!35!white] at ( 4.750,15.000) {\small\texttt{G190}};
\node at ( 5.875,15.000) {\small\texttt{G194}};
\node at ( 7.000,15.000) {\small\texttt{G195}};
\node at ( 8.125,15.000) {\small\texttt{G199}};
\node at ( 9.250,15.000) {\small\texttt{G200}};
\node at (10.375,15.000) {\small\texttt{G123}};
\node at (11.500,15.000) {\small\texttt{G124}};
\node at (12.625,15.000) {\small\texttt{G139}};
\node at ( 4.750,14.375) {\small\texttt{G203}};
\node at ( 5.875,14.375) {\small\texttt{G204}};
\node at ( 7.000,14.375) {\small\texttt{G205}};
\node at ( 8.125,14.375) {\small\texttt{G206}};
\node at ( 9.250,14.375) {\small\texttt{G207}};
\node[fill=yellow!35!white] at (10.375,14.375) {\small\texttt{~G93}};
\node at (11.500,14.375) {\small\texttt{~G95}};
\node at (12.625,14.375) {\small\texttt{G104}};
\node[fill=yellow!35!white] at ( 4.750,13.750) {\small\texttt{G114}};
\node at ( 5.875,13.750) {\small\texttt{G129}};
\node at ( 7.000,13.750) {\small\texttt{G135}};
\node[fill=yellow!35!white] at ( 8.125,13.750) {\small\texttt{G119}};
\node at ( 9.250,13.750) {\small\texttt{G130}};
\node[fill=yellow!35!white] at (10.375,13.750) {\small\texttt{~G80}};
\node[fill=yellow!35!white] at (11.500,13.750) {\small\texttt{~G97}};
\node[fill=yellow!35!white] at (12.625,13.750) {\small\texttt{G133}};
\node at ( 4.750,13.125) {\small\texttt{G145}};
\node at ( 5.875,13.125) {\small\texttt{G149}};
\node at ( 7.000,13.125) {\small\texttt{G162}};
\node at ( 8.125,13.125) {\small\texttt{G134}};
\node at ( 9.250,13.125) {\small\texttt{G142}};
\node at (10.375,13.125) {\small\texttt{~G81}};
\node at (11.500,13.125) {\small\texttt{G102}};
\node at (12.625,13.125) {\small\texttt{G153}};

\draw[thick,shift={(0,-1.25)}] 
(4.1875,14.0625) rectangle (13.1875,17.8125) 
( 7.5625,14.0625)--( 7.5625,15.3125) 
( 9.8125,14.0625)--( 9.8125,17.8125) (10.9375,14.0625)--(10.9375,15.3125) (12.0625,14.0625)--(12.0625,15.3125) 
( 4.1875,15.3125)--(13.1875,15.3125) 
( 9.8125,15.9375)--(13.1875,15.9375) 
( 4.1875,16.5625)--( 9.8125,16.5625)
( 4.1875,17.8125) rectangle (13.1875,19.6875)
( 8.6875,17.8125)--( 8.6875,18.4375)--(13.1875,18.4375)
(10.9375,17.8125)--(10.9375,18.4375)
( 9.8125,18.4375)--( 9.8125,19.6875)
( 9.8125,19.0625)--(13.1875,19.0625);
;

\node[fill=yellow!35!white] at (2.500,1.6875) {\small\texttt{G126}};
\node at (2.500,1.0625) {\small\texttt{G159}};
\node at (2.500,0.43750) {\small\texttt{G178}};

\node at (3.625,1.6875) {\small\texttt{G140}};
\node at (3.625,1.0625) {\small\texttt{G160}};
\node at (3.625,0.43750) {\small\texttt{G179}};

\node at (4.750,1.6875) {\small\texttt{G141}};
\node at (4.750,1.0625) {\small\texttt{G163}};
\node at (4.750,0.43750) {\small\texttt{G180}};

\node at (5.875,1.6875) {\small\texttt{G143}};
\node at (5.875,1.0625) {\small\texttt{G166}};
\node at (5.875,0.43750) {\small\texttt{G182}};

\node at (7.000,1.6875) {\small\texttt{G144}};
\node at (7.000,1.0625) {\small\texttt{G169}};
\node at (7.000,0.43750) {\small\texttt{G183}};

\node at (8.125,1.6875) {\small\texttt{G150}};
\node at (8.125,1.0625) {\small\texttt{G170}};
\node at (8.125,0.43750) {\small\texttt{G184}};

\node at (9.250,1.6875) {\small\texttt{G151}};
\node at (9.250,1.0625) {\small\texttt{G171}};
\node at (9.250,0.43750) {\small\texttt{G185}};

\node at (10.38,1.6875) {\small\texttt{G156}};
\node at (10.38,1.0625) {\small\texttt{G172}};
\node at (10.38,0.43750) {\small\texttt{G187}};

\node at (11.50,1.6875) {\small\texttt{G157}};
\node at (11.50,1.0625) {\small\texttt{G173}};
\node at (11.50,0.43750) {\small\texttt{G193}};

\node at (12.62,1.6875) {\small\texttt{G158}};
\node at (12.62,1.0625) {\small\texttt{G177}};

\draw[thick] (1.9375,2) rectangle (13.1875,.125);

\end{tikzpicture}

%% file: appendix.tex
\newpage

\appendix
\section{Connected graphs on six vertices realizing given ordered multiplicities}
For the following ordered multiplicities we list all connected graphs which can achieve a given multiplicity list.  Any graph which is underlined is a graph which has not yet been determined to be spectrally arbitrary for that ordered multiplicity list.  Any graph which is boxed is a graph which has been shown to \emph{not} be spectrally arbitrary for that multiplicity list.

\medskip\hrule\medskip
\noindent\textbf{111111}\medskip

\begin{tabular}{rrrrrrrrrrrr}
\texttt{G77} & \texttt{G78} & \texttt{G79} & \texttt{G80} & \texttt{G81} & \texttt{G83} & \texttt{G92} & \texttt{G93} & \texttt{G94} & \texttt{G95} & \texttt{G96} & \texttt{G97} \\
\texttt{G98} & \texttt{G99} & \texttt{G100} & \texttt{G102} & \texttt{G103} & \texttt{G104} & \texttt{G105} & \texttt{G111} & \texttt{G112} & \texttt{G113} & \texttt{G114} & \texttt{G115} \\
\texttt{G117} & \texttt{G118} & \texttt{G119} & \texttt{G120} & \texttt{G121} & \texttt{G122} & \texttt{G123} & \texttt{G124} & \texttt{G125} & \texttt{G126} & \texttt{G127} & \texttt{G128} \\
\texttt{G129} & \texttt{G130} & \texttt{G133} & \texttt{G134} & \texttt{G135} & \texttt{G136} & \texttt{G137} & \texttt{G138} & \texttt{G139} & \texttt{G140} & \texttt{G141} & \texttt{G142} \\
\texttt{G143} & \texttt{G144} & \texttt{G145} & \texttt{G146} & \texttt{G147} & \texttt{G148} & \texttt{G149} & \texttt{G150} & \texttt{G151} & \texttt{G152} & \texttt{G153} & \texttt{G154} \\
\texttt{G156} & \texttt{G157} & \texttt{G158} & \texttt{G159} & \texttt{G160} & \texttt{G161} & \texttt{G162} & \texttt{G163} & \texttt{G164} & \texttt{G165} & \texttt{G166} & \texttt{G167} \\
\texttt{G168} & \texttt{G169} & \texttt{G170} & \texttt{G171} & \texttt{G172} & \texttt{G173} & \texttt{G174} & \texttt{G175} & \texttt{G177} & \texttt{G178} & \texttt{G179} & \texttt{G180} \\
\texttt{G181} & \texttt{G182} & \texttt{G183} & \texttt{G184} & \texttt{G185} & \texttt{G186} & \texttt{G187} & \texttt{G188} & \texttt{G189} & \texttt{G190} & \texttt{G191} & \texttt{G192} \\
\texttt{G193} & \texttt{G194} & \texttt{G195} & \texttt{G196} & \texttt{G197} & \texttt{G198} & \texttt{G199} & \texttt{G200} & \texttt{G201} & \texttt{G202} & \texttt{G203} & \texttt{G204} \\
\texttt{G205} & \texttt{G206} & \texttt{G207} & \texttt{G208} 
\end{tabular}

\medskip\hrule\medskip
\noindent\textbf{11112 and 21111}\medskip

\begin{tabular}{rrrrrrrrrrrr}
\texttt{G92} & \texttt{G93} & \texttt{G94} & \texttt{G95} & \texttt{G96} & \texttt{G97} & \texttt{G98} & \texttt{G99} & \texttt{G100} & \texttt{G102} & \texttt{G103} & \texttt{G104} \\
\texttt{G105} & \texttt{G111} & \texttt{G112} & \texttt{G113} & \texttt{G114} & \texttt{G115} & \texttt{G117} & \texttt{G118} & \texttt{G119} & \texttt{G120} & \texttt{G121} & \texttt{G122} \\
\texttt{G123} & \texttt{G124} & \texttt{G125} & \texttt{G126} & \texttt{G127} & \texttt{G128} & \texttt{G129} & \texttt{G130} & \texttt{G133} & \texttt{G134} & \texttt{G135} & \texttt{G136} \\
\texttt{G137} & \texttt{G138} & \texttt{G139} & \texttt{G140} & \texttt{G141} & \texttt{G142} & \texttt{G143} & \texttt{G144} & \texttt{G145} & \texttt{G146} & \texttt{G147} & \texttt{G148} \\
\texttt{G149} & \texttt{G150} & \texttt{G151} & \texttt{G152} & \texttt{G153} & \texttt{G154} & \texttt{G156} & \texttt{G157} & \texttt{G158} & \texttt{G159} & \texttt{G160} & \texttt{G161} \\
\texttt{G162} & \texttt{G163} & \texttt{G164} & \texttt{G165} & \texttt{G166} & \texttt{G167} & \texttt{G168} & \texttt{G169} & \texttt{G170} & \texttt{G171} & \texttt{G172} & \texttt{G173} \\
\texttt{G174} & \texttt{G175} & \texttt{G177} & \texttt{G178} & \texttt{G179} & \texttt{G180} & \texttt{G181} & \texttt{G182} & \texttt{G183} & \texttt{G184} & \texttt{G185} & \texttt{G186} \\
\texttt{G187} & \texttt{G188} & \texttt{G189} & \texttt{G190} & \texttt{G191} & \texttt{G192} & \texttt{G193} & \texttt{G194} & \texttt{G195} & \texttt{G196} & \texttt{G197} & \texttt{G198} \\
\texttt{G199} & \texttt{G200} & \texttt{G201} & \texttt{G202} & \texttt{G203} & \texttt{G204} & \texttt{G205} & \texttt{G206} & \texttt{G207} & \texttt{G208} 
\end{tabular}

\medskip\hrule\medskip
\noindent\textbf{11121 and 12111}\medskip

\begin{tabular}{rrrrrrrrrrrr}
\texttt{G77} & \texttt{G78} & \texttt{G79} & \texttt{G80} & \texttt{G81} & \texttt{G92} & \texttt{G93} & \texttt{G94} & \texttt{G95} & \texttt{G96} & \texttt{G97} & \texttt{G98} \\
\texttt{G99} & \texttt{G100} & \texttt{G102} & \texttt{G103} & \texttt{G104} & \texttt{G105} & \texttt{G111} & \texttt{G112} & \texttt{G113} & \texttt{G114} & \texttt{G115} & \texttt{G117} \\
\texttt{G118} & \texttt{G119} & \texttt{G120} & \texttt{G121} & \texttt{G122} & \texttt{G123} & \texttt{G124} & \texttt{G125} & \texttt{G126} & \texttt{G127} & \texttt{G128} & \texttt{G129} \\
\texttt{G130} & \texttt{G133} & \texttt{G134} & \texttt{G135} & \texttt{G136} & \texttt{G137} & \texttt{G138} & \texttt{G139} & \texttt{G140} & \texttt{G141} & \texttt{G142} & \texttt{G143} \\
\texttt{G144} & \texttt{G145} & \texttt{G146} & \texttt{G147} & \texttt{G148} & \texttt{G149} & \texttt{G150} & \texttt{G151} & \texttt{G152} & \texttt{G153} & \texttt{G154} & \texttt{G156} \\
\texttt{G157} & \texttt{G158} & \texttt{G159} & \texttt{G160} & \texttt{G161} & \texttt{G162} & \texttt{G163} & \texttt{G164} & \texttt{G165} & \texttt{G166} & \texttt{G167} & \texttt{G168} \\
\texttt{G169} & \texttt{G170} & \texttt{G171} & \texttt{G172} & \texttt{G173} & \texttt{G174} & \texttt{G175} & \texttt{G177} & \texttt{G178} & \texttt{G179} & \texttt{G180} & \texttt{G181} \\
\texttt{G182} & \texttt{G183} & \texttt{G184} & \texttt{G185} & \texttt{G186} & \texttt{G187} & \texttt{G188} & \texttt{G189} & \texttt{G190} & \texttt{G191} & \texttt{G192} & \texttt{G193} \\
\texttt{G194} & \texttt{G195} & \texttt{G196} & \texttt{G197} & \texttt{G198} & \texttt{G199} & \texttt{G200} & \texttt{G201} & \texttt{G202} & \texttt{G203} & \texttt{G204} & \texttt{G205} \\
\texttt{G206} & \texttt{G207} & \texttt{G208} 
\end{tabular}

\medskip\hrule\medskip

\newpage

\noindent\textbf{11211}\medskip

\begin{tabular}{rrrrrrrrrrrr}
\texttt{G77} & \texttt{G78} & \texttt{G79} & \texttt{G80} & \texttt{G81} & \texttt{G92} & \texttt{G93} & \texttt{G94} & \texttt{G95} & \texttt{G96} & \texttt{G97} & \texttt{G98} \\
\texttt{G99} & \texttt{G100} & \texttt{G102} & \texttt{G103} & \texttt{G104} & \texttt{G105} & \texttt{G111} & \texttt{G112} & \texttt{G113} & \texttt{G114} & \texttt{G115} & \texttt{G117} \\
\texttt{G118} & \texttt{G119} & \texttt{G120} & \texttt{G121} & \texttt{G122} & \texttt{G123} & \texttt{G124} & \texttt{G125} & \texttt{G126} & \texttt{G127} & \texttt{G128} & \texttt{G129} \\
\texttt{G130} & \texttt{G133} & \texttt{G134} & \texttt{G135} & \texttt{G136} & \texttt{G137} & \texttt{G138} & \texttt{G139} & \texttt{G140} & \texttt{G141} & \texttt{G142} & \texttt{G143} \\
\texttt{G144} & \texttt{G145} & \texttt{G146} & \texttt{G147} & \texttt{G148} & \texttt{G149} & \texttt{G150} & \texttt{G151} & \texttt{G152} & \texttt{G153} & \texttt{G154} & \texttt{G156} \\
\texttt{G157} & \texttt{G158} & \texttt{G159} & \texttt{G160} & \texttt{G161} & \texttt{G162} & \texttt{G163} & \texttt{G164} & \texttt{G165} & \texttt{G166} & \texttt{G167} & \texttt{G168} \\
\texttt{G169} & \texttt{G170} & \texttt{G171} & \texttt{G172} & \texttt{G173} & \texttt{G174} & \texttt{G175} & \texttt{G177} & \texttt{G178} & \texttt{G179} & \texttt{G180} & \texttt{G181} \\
\texttt{G182} & \texttt{G183} & \texttt{G184} & \texttt{G185} & \texttt{G186} & \texttt{G187} & \texttt{G188} & \texttt{G189} & \texttt{G190} & \texttt{G191} & \texttt{G192} & \texttt{G193} \\
\texttt{G194} & \texttt{G195} & \texttt{G196} & \texttt{G197} & \texttt{G198} & \texttt{G199} & \texttt{G200} & \texttt{G201} & \texttt{G202} & \texttt{G203} & \texttt{G204} & \texttt{G205} \\
\texttt{G206} & \texttt{G207} & \texttt{G208} 
\end{tabular}

\medskip\hrule\medskip
\noindent\textbf{1122 and 2211}\medskip

\begin{tabular}{rrrrrrrrrrrr}
\texttt{G92} & \texttt{G93} & \texttt{G95} & \texttt{G96} & \texttt{G98} & \texttt{G99} & \texttt{G100} & \texttt{G103} & \texttt{G104} & \texttt{G105} & \texttt{G111} & \texttt{G112} \\
\texttt{G113} & \texttt{G114} & \texttt{G115} & \texttt{G117} & \texttt{G118} & \texttt{G119} & \texttt{G120} & \texttt{G121} & \texttt{G122} & \texttt{G123} & \texttt{G124} & \texttt{G125} \\
\texttt{G126} & \texttt{G127} & \texttt{G128} & \texttt{G129} & \texttt{G130} & \texttt{G133} & \texttt{G134} & \texttt{G135} & \texttt{G136} & \texttt{G137} & \texttt{G138} & \texttt{G139} \\
\texttt{G140} & \texttt{G141} & \texttt{G142} & \texttt{G143} & \texttt{G144} & \texttt{G145} & \texttt{G146} & \texttt{G147} & \texttt{G148} & \texttt{G149} & \texttt{G150} & \texttt{G151} \\
\texttt{G152} & \texttt{G153} & \texttt{G154} & \texttt{G156} & \texttt{G157} & \texttt{G158} & \texttt{G159} & \texttt{G160} & \texttt{G161} & \texttt{G162} & \texttt{G163} & \texttt{G164} \\
\texttt{G165} & \texttt{G166} & \texttt{G167} & \texttt{G168} & \texttt{G169} & \texttt{G170} & \texttt{G171} & \texttt{G172} & \texttt{G173} & \texttt{G174} & \texttt{G175} & \texttt{G177} \\
\texttt{G178} & \texttt{G179} & \texttt{G180} & \texttt{G181} & \texttt{G182} & \texttt{G183} & \texttt{G184} & \texttt{G185} & \texttt{G186} & \texttt{G187} & \texttt{G188} & \texttt{G189} \\
\texttt{G190} & \texttt{G191} & \texttt{G192} & \texttt{G193} & \texttt{G194} & \texttt{G195} & \texttt{G196} & \texttt{G197} & \texttt{G198} & \texttt{G199} & \texttt{G200} & \texttt{G201} \\
\texttt{G202} & \texttt{G203} & \texttt{G204} & \texttt{G205} & \texttt{G206} & \texttt{G207} & \texttt{G208} 
\end{tabular}

\medskip\hrule\medskip
\noindent\textbf{1221}\medskip

\begin{tabular}{rrrrrrrrrrrr}
\texttt{G79} & \texttt{G92} & \texttt{G93} & \texttt{G95} & \texttt{G96} & \texttt{G98} & \texttt{G99} & \underline{\texttt{G100}} & \texttt{G103} & \texttt{G104} & \texttt{G105} & \texttt{G111} \\
\texttt{G112} & \texttt{G113} & \texttt{G114} & \texttt{G115} & \texttt{G117} & \texttt{G118} & \texttt{G119} & \texttt{G120} & \texttt{G121} & \texttt{G122} & \texttt{G123} & \texttt{G124} \\
\texttt{G125} & \texttt{G126} & \texttt{G127} & \texttt{G128} & \texttt{G129} & \texttt{G130} & \texttt{G133} & \texttt{G134} & \texttt{G135} & \texttt{G136} & \texttt{G137} & \texttt{G138} \\
\texttt{G139} & \texttt{G140} & \texttt{G141} & \texttt{G142} & \texttt{G143} & \texttt{G144} & \texttt{G145} & \texttt{G146} & \texttt{G147} & \texttt{G148} & \texttt{G149} & \texttt{G150} \\
\texttt{G151} & \texttt{G152} & \texttt{G153} & \texttt{G154} & \texttt{G156} & \texttt{G157} & \texttt{G158} & \texttt{G159} & \texttt{G160} & \texttt{G161} & \texttt{G162} & \texttt{G163} \\
\texttt{G164} & \texttt{G165} & \texttt{G166} & \texttt{G167} & \texttt{G168} & \texttt{G169} & \texttt{G170} & \texttt{G171} & \texttt{G172} & \texttt{G173} & \texttt{G174} & \texttt{G175} \\
\texttt{G177} & \texttt{G178} & \texttt{G179} & \texttt{G180} & \texttt{G181} & \texttt{G182} & \texttt{G183} & \texttt{G184} & \texttt{G185} & \texttt{G186} & \texttt{G187} & \texttt{G188} \\
\texttt{G189} & \texttt{G190} & \texttt{G191} & \texttt{G192} & \texttt{G193} & \texttt{G194} & \texttt{G195} & \texttt{G196} & \texttt{G197} & \texttt{G198} & \texttt{G199} & \texttt{G200} \\
\texttt{G201} & \texttt{G202} & \texttt{G203} & \texttt{G204} & \texttt{G205} & \texttt{G206} & \texttt{G207} & \texttt{G208} 
\end{tabular}

\medskip\hrule\medskip
\noindent\textbf{1212 and 2121}\medskip

\begin{tabular}{rrrrrrrrrrrr}
\texttt{G92} & \texttt{G93} & \underline{\texttt{G94}} & \texttt{G95} & \texttt{G96} & \texttt{G98} & \texttt{G99} & \texttt{G100} & \texttt{G103} & \texttt{G104} & \texttt{G111} &
\texttt{G112} \\ \texttt{G113} & \texttt{G114} & \texttt{G115} & \texttt{G117} & \texttt{G118} & \texttt{G119} & \texttt{G120} & \texttt{G121} & \texttt{G122} & \texttt{G123} & \texttt{G124} &
\texttt{G125} \\ \texttt{G126} & \texttt{G127} & \texttt{G128} & \texttt{G129} & \texttt{G130} & \texttt{G133} & \texttt{G134} & \texttt{G135} & \texttt{G136} & \texttt{G137} & \texttt{G138} &
\texttt{G139} \\ \texttt{G140} & \texttt{G141} & \texttt{G142} & \texttt{G143} & \texttt{G144} & \texttt{G145} & \texttt{G146} & \texttt{G147} & \texttt{G148} & \texttt{G149} & \texttt{G150} &
\texttt{G151} \\ \texttt{G152} & \texttt{G153} & \texttt{G154} & \texttt{G156} & \texttt{G157} & \texttt{G158} & \texttt{G159} & \texttt{G160} & \texttt{G161} & \texttt{G162} & \texttt{G163} &
\texttt{G164} \\ \texttt{G165} & \texttt{G166} & \texttt{G167} & \texttt{G168} & \texttt{G169} & \texttt{G170} & \texttt{G171} & \texttt{G172} & \texttt{G173} & \texttt{G174} & \texttt{G175} &
\texttt{G177} \\ \texttt{G178} & \texttt{G179} & \texttt{G180} & \texttt{G181} & \texttt{G182} & \texttt{G183} & \texttt{G184} & \texttt{G185} & \texttt{G186} & \texttt{G187} & \texttt{G188} &
\texttt{G189} \\ \texttt{G190} & \texttt{G191} & \texttt{G192} & \texttt{G193} & \texttt{G194} & \texttt{G195} & \texttt{G196} & \texttt{G197} & \texttt{G198} & \texttt{G199} & \texttt{G200} &
\texttt{G201} \\ \texttt{G202} & \texttt{G203} & \texttt{G204} & \texttt{G205} & \texttt{G206} & \texttt{G207} & \texttt{G208} 
\end{tabular}

\medskip\hrule\medskip

\newpage

\noindent\textbf{2112}\medskip

\begin{tabular}{rrrrrrrrrrrr}
\texttt{G96} & \texttt{G98} & \texttt{G99} & \texttt{G103} & \texttt{G105} & \texttt{G111} & \texttt{G112} & \texttt{G113} & \texttt{G114} & \texttt{G115} & \texttt{G117} & \texttt{G118} \\
\texttt{G119} & \texttt{G120} & \texttt{G121} & \texttt{G122} & \texttt{G123} & \texttt{G124} & \texttt{G125} & \texttt{G126} & \underline{\texttt{G127}} & \texttt{G128} & \texttt{G129} & \texttt{G130} \\
\texttt{G133} & \texttt{G134} & \texttt{G135} & \texttt{G136} & \texttt{G137} & \texttt{G138} & \texttt{G139} & \texttt{G140} & \texttt{G141} & \texttt{G142} & \texttt{G143} & \texttt{G144} \\
\texttt{G145} & \texttt{G146} & \texttt{G147} & \texttt{G148} & \texttt{G149} & \texttt{G150} & \texttt{G151} & \texttt{G152} & \texttt{G153} & \texttt{G154} & \texttt{G156} & \texttt{G157} \\
\texttt{G158} & \texttt{G159} & \texttt{G160} & \texttt{G161} & \texttt{G162} & \texttt{G163} & \texttt{G164} & \texttt{G165} & \texttt{G166} & \texttt{G167} & \texttt{G168} & \texttt{G169} \\
\texttt{G170} & \texttt{G171} & \texttt{G172} & \texttt{G173} & \texttt{G174} & \texttt{G175} & \texttt{G177} & \texttt{G178} & \texttt{G179} & \texttt{G180} & \texttt{G181} & \texttt{G182} \\
\texttt{G183} & \texttt{G184} & \texttt{G185} & \texttt{G186} & \texttt{G187} & \texttt{G188} & \texttt{G189} & \texttt{G190} & \texttt{G191} & \texttt{G192} & \texttt{G193} & \texttt{G194} \\
\texttt{G195} & \texttt{G196} & \texttt{G197} & \texttt{G198} & \texttt{G199} & \texttt{G200} & \texttt{G201} & \texttt{G202} & \texttt{G203} & \texttt{G204} & \texttt{G205} & \texttt{G206} \\
\texttt{G207} & \texttt{G208} 
\end{tabular}

\medskip\hrule\medskip
\noindent\textbf{222}\medskip

\begin{tabular}{rrrrrrrrrrrr}
\texttt{G96} & \texttt{G99} & \texttt{G105} & \texttt{G111} & \texttt{G114} & \underline{\texttt{G115}} & \texttt{G117} & \texttt{G118} & \texttt{G121} & \underline{\texttt{G125}} & \texttt{G126} & \texttt{G127} \\
\texttt{G128} & \underline{\texttt{G129}} & \texttt{G133} & \texttt{G135} & \texttt{G136} & \texttt{G137} & \underline{\texttt{G138}} & \texttt{G140} & \texttt{G141} & \texttt{G143} & \texttt{G144} & \texttt{G145} \\
\texttt{G146} & \texttt{G147} & \texttt{G148} & \texttt{G149} & \texttt{G150} & \texttt{G151} & \texttt{G152} & \texttt{G153} & \texttt{G154} & \texttt{G156} & \texttt{G157} & \texttt{G158} \\
\texttt{G159} & \texttt{G160} & \texttt{G161} & \texttt{G162} & \texttt{G163} & \texttt{G164} & \texttt{G165} & \texttt{G166} & \texttt{G167} & \texttt{G168} & \texttt{G169} & \texttt{G170} \\
\texttt{G171} & \texttt{G172} & \texttt{G173} & \texttt{G174} & \texttt{G175} & \texttt{G177} & \texttt{G178} & \texttt{G179} & \texttt{G180} & \texttt{G181} & \texttt{G182} & \texttt{G183} \\
\texttt{G184} & \texttt{G185} & \texttt{G186} & \texttt{G187} & \texttt{G188} & \texttt{G189} & \texttt{G190} & \texttt{G191} & \texttt{G192} & \texttt{G193} & \texttt{G194} & \texttt{G195} \\
\texttt{G196} & \texttt{G197} & \texttt{G198} & \texttt{G199} & \texttt{G200} & \texttt{G201} & \texttt{G202} & \texttt{G203} & \texttt{G204} & \texttt{G205} & \texttt{G206} & \texttt{G207} \\
\texttt{G208} 
\end{tabular}

\medskip\hrule\medskip
\noindent\textbf{1113 and 3111}\medskip

\begin{tabular}{rrrrrrrrrrrr}
\texttt{G117} & \texttt{G119} & \texttt{G126} & \texttt{G130} & \texttt{G133} & \texttt{G134} & \texttt{G140} & \texttt{G141} & \texttt{G142} & \texttt{G143} & \texttt{G144} & \texttt{G150} \\
\texttt{G151} & \texttt{G153} & \texttt{G154} & \texttt{G156} & \texttt{G157} & \texttt{G158} & \texttt{G159} & \texttt{G160} & \texttt{G163} & \texttt{G165} & \texttt{G166} & \texttt{G168} \\
\texttt{G169} & \texttt{G170} & \texttt{G171} & \texttt{G172} & \texttt{G173} & \texttt{G174} & \texttt{G175} & \texttt{G177} & \texttt{G178} & \texttt{G179} & \texttt{G180} & \texttt{G181} \\
\texttt{G182} & \texttt{G183} & \texttt{G184} & \texttt{G185} & \texttt{G186} & \texttt{G187} & \texttt{G188} & \texttt{G189} & \texttt{G190} & \texttt{G191} & \texttt{G192} & \texttt{G193} \\
\texttt{G194} & \texttt{G195} & \texttt{G196} & \texttt{G197} & \texttt{G198} & \texttt{G199} & \texttt{G200} & \texttt{G201} & \texttt{G202} & \texttt{G203} & \texttt{G204} & \texttt{G205} \\
\texttt{G206} & \texttt{G207} & \texttt{G208} 
\end{tabular}

\medskip\hrule\medskip
\noindent\textbf{1131 and 1311}\medskip

\begin{tabular}{rrrrrrrrrrrr}
\texttt{G77} & \texttt{G78} & \underline{\texttt{G92}} & \texttt{G100} & \texttt{G114} & \underline{\texttt{G117}} & \texttt{G119} & \texttt{G121} & \texttt{G125} & \texttt{G126} & \underline{\texttt{G129}} & \underline{\texttt{G130}} \\
\texttt{G133} & \texttt{G134} & \texttt{G135} & \texttt{G138} & \texttt{G140} & \texttt{G141} & \texttt{G142} & \texttt{G143} & \texttt{G144} & \underline{\texttt{G145}} & \texttt{G146} & \texttt{G149} \\
\underline{\texttt{G150}} & \underline{\texttt{G151}} & \underline{\texttt{G153}} & \texttt{G154} & \texttt{G156} & \texttt{G157} & \texttt{G158} & \texttt{G159} & \texttt{G160} & \texttt{G161} & \texttt{G162} & \underline{\texttt{G163}} \\
\texttt{G165} & \texttt{G166} & \texttt{G168} & \texttt{G169} & \texttt{G170} & \underline{\texttt{G171}} & \texttt{G172} & \texttt{G173} & \underline{\texttt{G174}} & \texttt{G175} & \texttt{G177} & \texttt{G178} \\
\texttt{G179} & \texttt{G180} & \texttt{G181} & \texttt{G182} & \texttt{G183} & \texttt{G184} & \texttt{G185} & \texttt{G186} & \underline{\texttt{G187}} & \texttt{G188} & \texttt{G189} & \texttt{G190} \\
\texttt{G191} & \texttt{G192} & \texttt{G193} & \texttt{G194} & \texttt{G195} & \texttt{G196} & \texttt{G197} & \texttt{G198} & \texttt{G199} & \texttt{G200} & \texttt{G201} & \texttt{G202} \\
\texttt{G203} & \texttt{G204} & \texttt{G205} & \texttt{G206} & \texttt{G207} & \texttt{G208} 
\end{tabular}

\medskip\hrule\medskip
\noindent\textbf{123 and 321}\medskip

\begin{tabular}{rrrrrrrrrrrr}
\texttt{G117} & \texttt{G126} & \texttt{G133} & \texttt{G140} & \texttt{G141} & \texttt{G143} & \texttt{G144} & \texttt{G150} & \underline{\texttt{G151}} & \texttt{G153} & \texttt{G154} & \texttt{G156} \\
\texttt{G157} & \texttt{G158} & \texttt{G159} & \texttt{G160} & \underline{\texttt{G163}} & \texttt{G165} & \texttt{G166} & \texttt{G168} & \texttt{G169} & \texttt{G170} & \underline{\texttt{G171}} & \texttt{G172} \\
\texttt{G173} & \texttt{G174} & \texttt{G175} & \texttt{G177} & \texttt{G178} & \texttt{G179} & \texttt{G180} & \texttt{G181} & \texttt{G182} & \texttt{G183} & \texttt{G184} & \texttt{G185} \\
\texttt{G186} & \underline{\texttt{G187}} & \texttt{G188} & \texttt{G189} & \texttt{G190} & \texttt{G191} & \texttt{G192} & \texttt{G193} & \texttt{G194} & \texttt{G195} & \texttt{G196} & \texttt{G197} \\
\texttt{G198} & \texttt{G199} & \texttt{G200} & \texttt{G201} & \texttt{G202} & \texttt{G203} & \texttt{G204} & \texttt{G205} & \texttt{G206} & \texttt{G207} & \texttt{G208} 
\end{tabular}

\medskip\hrule\medskip

\newpage

\noindent\textbf{132 and 231}\medskip

\begin{tabular}{rrrrrrrrrrrr}
\texttt{G92} & \underline{\texttt{G114}} & \texttt{G126} & \underline{\texttt{G129}} & \texttt{G133} & \underline{\texttt{G135}} & \texttt{G140} & \texttt{G141} & \texttt{G143} & \texttt{G144} & \underline{\texttt{G145}} & \underline{\texttt{G146}} \\
\underline{\texttt{G149}} & \underline{\texttt{G150}} & \underline{\texttt{G151}} & \underline{\texttt{G153}} & \underline{\texttt{G154}} & \texttt{G156} & \texttt{G157} & \texttt{G158} & \texttt{G159} & \texttt{G160} & \texttt{G161} & \underline{\texttt{G162}} \\
\underline{\texttt{G163}} & \texttt{G165} & \texttt{G166} & \texttt{G168} & \underline{\texttt{G169}} & \texttt{G170} & \underline{\texttt{G171}} & \texttt{G172} & \texttt{G173} & \underline{\texttt{G174}} & \underline{\texttt{G175}} & \texttt{G177} \\
\texttt{G178} & \texttt{G179} & \texttt{G180} & \texttt{G181} & \texttt{G182} & \texttt{G183} & \texttt{G184} & \texttt{G185} & \texttt{G186} & \underline{\texttt{G187}} & \texttt{G188} & \texttt{G189} \\
\texttt{G190} & \texttt{G191} & \texttt{G192} & \texttt{G193} & \texttt{G194} & \texttt{G195} & \texttt{G196} & \texttt{G197} & \texttt{G198} & \texttt{G199} & \texttt{G200} & \texttt{G201} \\
\texttt{G202} & \texttt{G203} & \texttt{G204} & \texttt{G205} & \texttt{G206} & \texttt{G207} & \texttt{G208} 
\end{tabular}

\medskip\hrule\medskip
\noindent\textbf{213 and 312}\medskip

\begin{tabular}{rrrrrrrrrrrr}
\texttt{G126} & \texttt{G140} & \texttt{G141} & \texttt{G143} & \texttt{G144} & \texttt{G150} & \texttt{G151} & \texttt{G154} & \texttt{G156} & \texttt{G157} & \texttt{G158} & \texttt{G159} \\
\texttt{G160} & \texttt{G163} & \texttt{G165} & \texttt{G166} & \texttt{G168} & \texttt{G169} & \texttt{G170} & \texttt{G171} & \texttt{G172} & \texttt{G173} & \texttt{G174} & \texttt{G175} \\
\texttt{G177} & \texttt{G178} & \texttt{G179} & \texttt{G180} & \texttt{G181} & \texttt{G182} & \texttt{G183} & \texttt{G184} & \texttt{G185} & \texttt{G186} & \texttt{G187} & \texttt{G188} \\
\texttt{G189} & \texttt{G190} & \texttt{G191} & \texttt{G192} & \texttt{G193} & \texttt{G194} & \texttt{G195} & \texttt{G196} & \texttt{G197} & \texttt{G198} & \texttt{G199} & \texttt{G200} \\
\texttt{G201} & \texttt{G202} & \texttt{G203} & \texttt{G204} & \texttt{G205} & \texttt{G206} & \texttt{G207} & \texttt{G208} 
\end{tabular}

\medskip\hrule\medskip
\noindent\textbf{33}\medskip

\begin{tabular}{rrrrrrrrrrrr}
\texttt{G154} & \texttt{G168} & \texttt{G174} & \texttt{G175} & \texttt{G181} & \texttt{G186} & \texttt{G188} & \texttt{G190} & \texttt{G192} & \texttt{G194} & \texttt{G195} & \texttt{G196} \\
\texttt{G197} & \texttt{G198} & \texttt{G199} & \texttt{G200} & \texttt{G201} & \texttt{G202} & \texttt{G203} & \texttt{G204} & \texttt{G205} & \texttt{G206} & \texttt{G207} & \texttt{G208} \\
\end{tabular}

\medskip\hrule\medskip
\noindent\textbf{114 and 411}\medskip

\begin{tabular}{rrrrrrrrrrrr}
\texttt{G165} & \texttt{G190} & \texttt{G191} & \texttt{G194} & \texttt{G195} & \texttt{G199} & \texttt{G200} & \texttt{G203} & \texttt{G204} & \texttt{G205} & \texttt{G206} & \texttt{G207} \\
\texttt{G208} 
\end{tabular}

\medskip\hrule\medskip
\noindent\textbf{141}\medskip

\begin{tabular}{rrrrrrrrrrrr}
\texttt{G77} & \texttt{G146} & \texttt{G161} & \texttt{G165} & \boxed{\texttt{G175}} & \texttt{G189} & \texttt{G190} & \texttt{G191} & \texttt{G194} & \texttt{G195} & \texttt{G197} & \texttt{G199} \\
\texttt{G200} & \texttt{G201} & \texttt{G203} & \texttt{G204} & \texttt{G205} & \texttt{G206} & \texttt{G207} & \texttt{G208} 
\end{tabular}

\medskip\hrule\medskip
\noindent\textbf{24 and 42}\medskip

\begin{tabular}{rrrrrrrrrrrr}
\texttt{G190} & \texttt{G194} & \texttt{G195} & \texttt{G199} & \texttt{G200} & \texttt{G203} & \texttt{G204} & \texttt{G205} & \texttt{G206} & \texttt{G207} & \texttt{G208} 
\end{tabular}

\medskip\hrule\medskip
\noindent\textbf{15 and 51}\medskip

\begin{tabular}{rrrrrrrrrrrr}
\texttt{G208} 
\end{tabular}

\newpage

\section{Attainable ordered multiplicitie lists}

For each connected graph on six or fewer vertices we give the Atlas of Graph numbering, draw the graph, and list all attainable ordered multiplicity lists.

\GG{1}{1}

\GG{3}{11}

\GG{6}{111}

\GG{7}{111, 12, 21}

\GG{13}{1111, 121}

\GG{14}{1111}

\GG{15}{1111, 121, 112, 211}

\GG{16}{1111, 121, 112, 211, 22}

\GG{17}{1111, 121, 112, 211, 22}

\GG{18}{1111, 121, 112, 211, 22, 13, 31}

\GG{29}{11111, 1121, 1211, 131}

\GG{30}{11111, 1121, 1211}

\GG{31}{11111}

\GG{34}{11111, 1121, 1211, 1112, 2111, 122, 221}

\GG{35}{11111, 1121, 1211, 1112, 2111}

\GG{36}{11111, 1121, 1211, 1112, 2111}

\GG{37}{11111, 1121, 1211, 1112, 2111, 122, 221, 212}

\GG{38}{11111, 1121, 1211, 1112, 2111, 122, 221}

\GG{40}{11111, 1121, 1211, 1112, 2111, 122, 221, 212}

\GG{41}{11111, 1121, 1211, 1112, 2111, 122, 221, 212}

\GG{42}{11111, 1121, 1211, 1112, 2111, 122, 221, 212, 113, 131, 311}

\GG{43}{11111, 1121, 1211, 1112, 2111, 122, 221, 212}

\GG{44}{11111, 1121, 1211, 1112, 2111, 122, 221, 212, 131}

\GG{45}{11111, 1121, 1211, 1112, 2111, 122, 221, 212, 131, 113, 311}

\GG{46}{11111, 1121, 1211, 1112, 2111, 122, 221, 212, 131}

\GG{47}{11111, 1121, 1211, 1112, 2111, 122, 221, 212}

\GG{48}{11111, 1121, 1211, 1112, 2111, 122, 221, 212, 131, 113, 311, 23, 32}

\GG{49}{11111, 1121, 1211, 1112, 2111, 122, 221, 212, 131, 113, 311, 23, 32}

\GG{50}{11111, 1121, 1211, 1112, 2111, 122, 221, 212, 131, 113, 311, 23, 32}

\GG{51}{11111, 1121, 1211, 1112, 2111, 122, 221, 212, 131, 113, 311, 23, 32}

\GG{52}{11111, 1121, 1211, 1112, 2111, 122, 221, 212, 131, 113, 311, 23, 32, 14, 41}

\GG{77}{111111, 11121, 12111, 11211, 1131, 1311, 141}

\GG{78}{111111, 11121, 12111, 11211, 1131, 1311}

\GG{79}{111111, 11121, 12111, 11211, 1221}

\GG{80}{111111, 11121, 12111, 11211}

\GG{81}{111111, 11121, 12111, 11211}

\GG{83}{111111}

\GG{92}{111111, 11112, 21111, 11121, 12111, 11211, 1122, 2211, 1212, 2121, 1221, \underline{1131}, \underline{1311}, 132, 231}

\GG{93}{111111, 11112, 21111, 11121, 12111, 11211, 1122, 2211, 1212, 2121, 1221}

\GG{94}{111111, 11112, 21111, 11121, 12111, 11211, \underline{1212}, \underline{2121}}

\GG{95}{111111, 11112, 21111, 11121, 12111, 11211, 1122, 2211, 1212, 2121, 1221}

\GG{96}{111111, 11112, 21111, 11121, 12111, 11211, 1122, 2211, 1212, 2121, 1221, 2112, 222}

\GG{97}{111111, 11112, 21111, 11121, 12111, 11211}

\GG{98}{111111, 11112, 21111, 11121, 12111, 11211, 1122, 2211, 1212, 2121, 1221, 2112}

\GG{99}{111111, 11112, 21111, 11121, 12111, 11211, 1122, 2211, 1212, 2121, 1221, 2112, 222}

\GG{100}{111111, 11112, 21111, 11121, 12111, 11211, 1122, 2211, 1212, 2121, \underline{1221}, 1131, 1311}

\GG{102}{111111, 11112, 21111, 11121, 12111, 11211}

\GG{103}{111111, 11112, 21111, 11121, 12111, 11211, 1122, 2211, 1212, 2121, 1221, 2112}

\GG{104}{111111, 11112, 21111, 11121, 12111, 11211, 1122, 2211, 1212, 2121, 1221}

\GG{105}{111111, {11112}, {21111}, {11121}, {12111}, {11211}, {1122}, {2211}, {1221}, {2112}, 222}

\GG{111}{111111, 11112, 21111, 11121, 12111, 11211, 1122, 2211, 1212, 2121, 1221, 2112, 222}

\GG{112}{111111, 11112, 21111, 11121, 12111, 11211, 1122, 2211, 1212, 2121, 1221, 2112}

\GG{113}{111111, 11112, 21111, 11121, 12111, 11211, 1122, 2211, 1212, 2121, 1221, 2112}

\GG{114}{111111, 11112, 21111, 11121, 12111, 11211, 1122, 2211, 1212, 2121, 1221, 2112, 222, 1131, 1311, \underline{132}, \underline{231}}

\GG{115}{111111, 11112, 21111, 11121, 12111, 11211, 1122, 2211, 1212, 2121, 1221, 2112, \underline{222}}

\GG{117}{111111, 11112, 21111, 11121, 12111, 11211, 1122, 2211, 1212, 2121, 1221, 2112, 222, 1113, 3111, \underline{1131}, \underline{1311}, 123, 321}

\GG{118}{111111, 11112, 21111, 11121, 12111, 11211, 1122, 2211, 1212, 2121, 1221, 2112, 222}

\GG{119}{111111, 11112, 21111, 11121, 12111, 11211, 1122, 2211, 1212, 2121, 1221, 2112, 1113, 3111, 1131, 1311}

\GG{120}{111111, 11112, 21111, 11121, 12111, 11211, 1122, 2211, 1212, 2121, 1221, 2112}

\GG{121}{111111, 11112, 21111, 11121, 12111, 11211, 1122, 2211, 1212, 2121, 1221, 2112, 222, 1131, 1311}

\GG{122}{111111, 11112, 21111, 11121, 12111, 11211, 1122, 2211, 1212, 2121, 1221, 2112}

\GG{123}{111111, 11112, 21111, 11121, 12111, 11211, 1122, 2211, 1212, 2121, 1221, 2112}

\GG{124}{111111, 11112, 21111, 11121, 12111, 11211, 1122, 2211, 1212, 2121, 1221, 2112}

\GG{125}{111111, 11112, 21111, 11121, 12111, 11211, 1122, 2211, 1212, 2121, 1221, 2112, \underline{222}, 1131, 1311}

\GG{126}{111111, 11112, 21111, 11121, 12111, 11211, 1122, 2211, 1212, 2121, 1221, 2112, 222, 1113, 3111, 1131, 1311, 123, 321, 132, 231, 213, 312}

\GG{127}{111111, 11112, 21111, 11121, 12111, 11211, 1122, 2211, 1212, 2121, 1221, \underline{2112}, 222}

\GG{128}{111111, 11112, 21111, 11121, 12111, 11211, 1122, 2211, 1212, 2121, 1221, 2112, 222}

\GG{129}{111111, 11112, 21111, 11121, 12111, 11211, 1122, 2211, 1212, 2121, 1221, 2112, \underline{222}, \underline{1131}, \underline{1311}, \underline{132}, \underline{231}}

\GG{130}{111111, 11112, 21111, 11121, 12111, 11211, 1122, 2211, 1212, 2121, 1221, 2112, 1113, 3111, \underline{1131}, \underline{1311}}

\GG{133}{111111, 11112, 21111, 11121, 12111, 11211, 1122, 2211, 1212, 2121, 1221, 2112, 222, 1113, 3111, 1131, 1311, 123, 321, 132, 231}

\GG{134}{111111, 11112, 21111, 11121, 12111, 11211, 1122, 2211, 1212, 2121, 1221, 2112, 1113, 3111, 1131, 1311}

\GG{135}{111111, 11112, 21111, 11121, 12111, 11211, 1122, 2211, 1212, 2121, 1221, 2112, 222, 1131, 1311, \underline{132}, \underline{231}}

\GG{136}{111111, 11112, 21111, 11121, 12111, 11211, 1122, 2211, 1212, 2121, 1221, 2112, 222}

\GG{137}{111111, 11112, 21111, 11121, 12111, 11211, 1122, 2211, 1212, 2121, 1221, 2112, 222}

\GG{138}{111111, 11112, 21111, 11121, 12111, 11211, 1122, 2211, 1212, 2121, 1221, 2112, \underline{222}, 1131, 1311}

\GG{139}{111111, 11112, 21111, 11121, 12111, 11211, 1122, 2211, 1212, 2121, 1221, 2112}

\GG{140}{111111, 11112, 21111, 11121, 12111, 11211, 1122, 2211, 1212, 2121, 1221, 2112, 222, 1113, 3111, 1131, 1311, 123, 321, 132, 231, 213, 312}

\GG{141}{111111, 11112, 21111, 11121, 12111, 11211, 1122, 2211, 1212, 2121, 1221, 2112, 222, 1113, 3111, 1131, 1311, 123, 321, 132, 231, 213, 312}

\GG{142}{111111, 11112, 21111, 11121, 12111, 11211, 1122, 2211, 1212, 2121, 1221, 2112, 1113, 3111, 1131, 1311}

\GG{143}{111111, 11112, 21111, 11121, 12111, 11211, 1122, 2211, 1212, 2121, 1221, 2112, 222, 1113, 3111, 1131, 1311, 123, 321, 132, 231, 213, 312}

\GG{144}{111111, 11112, 21111, 11121, 12111, 11211, 1122, 2211, 1212, 2121, 1221, 2112, 222, 1113, 3111, 1131, 1311, 123, 321, 132, 231, 213, 312}

\GG{145}{111111, 11112, 21111, 11121, 12111, 11211, 1122, 2211, 1212, 2121, 1221, 2112, 222, \underline{1131}, \underline{1311}, \underline{132}, \underline{231}}

\GG{146}{111111, 11112, 21111, 11121, 12111, 11211, 1122, 2211, 1212, 2121, 1221, 2112, 222, 1131, 1311, \underline{132}, \underline{231}, 141}

\GG{147}{111111, 11112, 21111, 11121, 12111, 11211, 1122, 2211, 1212, 2121, 1221, 2112, 222}

\GG{148}{111111, 11112, 21111, 11121, 12111, 11211, 1122, 2211, 1212, 2121, 1221, 2112, 222}

\GG{149}{111111, 11112, 21111, 11121, 12111, 11211, 1122, 2211, 1212, 2121, 1221, 2112, 222, 1131, 1311, \underline{132}, \underline{231}}

\GG{150}{111111, 11112, 21111, 11121, 12111, 11211, 1122, 2211, 1212, 2121, 1221, 2112, 222, 1113, 3111, \underline{1131}, \underline{1311}, 123, 321, \underline{132}, \underline{231}, 213, 312}

\GG{151}{111111, 11112, 21111, 11121, 12111, 11211, 1122, 2211, 1212, 2121, 1221, 2112, 222, 1113, 3111, \underline{1131}, \underline{1311}, \underline{123}, \underline{321}, \underline{132}, \underline{231}, 213, 312}

\GG{152}{111111, 11112, 21111, 11121, 12111, 11211, 1122, 2211, 1212, 2121, 1221, 2112, 222}

\GG{153}{111111, 11112, 21111, 11121, 12111, 11211, 1122, 2211, 1212, 2121, 1221, 2112, 222, 1113, 3111, \underline{1131}, \underline{1311}, 123, 321, \underline{132}, \underline{231}}

\GG{154}{111111, 11112, 21111, 11121, 12111, 11211, 1122, 2211, 1212, 2121, 1221, 2112, 222, 1113, 3111, 1131, 1311, 123, 321, \underline{132}, \underline{231}, 213, 312, 33}

\GG{156}{111111, 11112, 21111, 11121, 12111, 11211, 1122, 2211, 1212, 2121, 1221, 2112, 222, 1113, 3111, 1131, 1311, 123, 321, 132, 231, 213, 312}

\GG{157}{111111, 11112, 21111, 11121, 12111, 11211, 1122, 2211, 1212, 2121, 1221, 2112, 222, 1113, 3111, 1131, 1311, 123, 321, 132, 231, 213, 312}

\GG{158}{111111, 11112, 21111, 11121, 12111, 11211, 1122, 2211, 1212, 2121, 1221, 2112, 222, 1113, 3111, 1131, 1311, 123, 321, 132, 231, 213, 312}

\GG{159}{111111, 11112, 21111, 11121, 12111, 11211, 1122, 2211, 1212, 2121, 1221, 2112, 222, 1113, 3111, 1131, 1311, 123, 321, 132, 231, 213, 312}

\GG{160}{111111, 11112, 21111, 11121, 12111, 11211, 1122, 2211, 1212, 2121, 1221, 2112, 222, 1113, 3111, 1131, 1311, 123, 321, 132, 231, 213, 312}

\GG{161}{111111, 11112, 21111, 11121, 12111, 11211, 1122, 2211, 1212, 2121, 1221, 2112, 222, 1131, 1311, 132, 231, 141}

\GG{162}{111111, 11112, 21111, 11121, 12111, 11211, 1122, 2211, 1212, 2121, 1221, 2112, 222, 1131, 1311, \underline{132}, \underline{231}}

\GG{163}{111111, 11112, 21111, 11121, 12111, 11211, 1122, 2211, 1212, 2121, 1221, 2112, 222, 1113, 3111, \underline{1131}, \underline{1311}, \underline{123}, \underline{321}, \underline{132}, \underline{231}, 213, 312}

\GG{164}{111111, 11112, 21111, 11121, 12111, 11211, 1122, 2211, 1212, 2121, 1221, 2112, 222}

\GG{165}{111111, 11112, 21111, 11121, 12111, 11211, 1122, 2211, 1212, 2121, 1221, 2112, 222, 1113, 3111, 1131, 1311, 123, 321, 132, 231, 213, 312, 114, 411, 141}

\GG{166}{111111, 11112, 21111, 11121, 12111, 11211, 1122, 2211, 1212, 2121, 1221, 2112, 222, 1113, 3111, 1131, 1311, 123, 321, 132, 231, 213, 312}

\GG{167}{111111, 11112, 21111, 11121, 12111, 11211, 1122, 2211, 1212, 2121, 1221, 2112, 222}

\GG{168}{111111, 11112, 21111, 11121, 12111, 11211, 1122, 2211, 1212, 2121, 1221, 2112, 222, 1113, 3111, 1131, 1311, 123, 321, 132, 231, 213, 312, 33}

\GG{169}{111111, 11112, 21111, 11121, 12111, 11211, 1122, 2211, 1212, 2121, 1221, 2112, 222, 1113, 3111, 1131, 1311, 123, 321, \underline{132}, \underline{231}, 213, 312}

\GG{170}{111111, 11112, 21111, 11121, 12111, 11211, 1122, 2211, 1212, 2121, 1221, 2112, 222, 1113, 3111, 1131, 1311, 123, 321, 132, 231, 213, 312}

\GG{171}{111111, 11112, 21111, 11121, 12111, 11211, 1122, 2211, 1212, 2121, 1221, 2112, 222, 1113, 3111, \underline{1131}, \underline{1311}, \underline{123}, \underline{321}, \underline{132}, \underline{231}, 213, 312}

\GG{172}{111111, 11112, 21111, 11121, 12111, 11211, 1122, 2211, 1212, 2121, 1221, 2112, 222, 1113, 3111, 1131, 1311, 123, 321, 132, 231, 213, 312}

\GG{173}{111111, 11112, 21111, 11121, 12111, 11211, 1122, 2211, 1212, 2121, 1221, 2112, 222, 1113, 3111, 1131, 1311, 123, 321, 132, 231, 213, 312}

\GG{174}{111111, 11112, 21111, 11121, 12111, 11211, 1122, 2211, 1212, 2121, 1221, 2112, 222, 1113, 3111, \underline{1131}, \underline{1311}, 123, 321, \underline{132}, \underline{231}, 213, 312, 33}

\GG{175}{111111, 11112, 21111, 11121, 12111, 11211, 1122, 2211, 1212, 2121, 1221, 2112, 222, 1113, 3111, 1131, 1311, 123, 321, \underline{132}, \underline{231}, 213, 312, 33, \boxed{141}}

\GG{177}{111111, 11112, 21111, 11121, 12111, 11211, 1122, 2211, 1212, 2121, 1221, 2112, 222, 1113, 3111, 1131, 1311, 123, 321, 132, 231, 213, 312}

\GG{178}{111111, 11112, 21111, 11121, 12111, 11211, 1122, 2211, 1212, 2121, 1221, 2112, 222, 1113, 3111, 1131, 1311, 123, 321, 132, 231, 213, 312}

\GG{179}{111111, 11112, 21111, 11121, 12111, 11211, 1122, 2211, 1212, 2121, 1221, 2112, 222, 1113, 3111, 1131, 1311, 123, 321, 132, 231, 213, 312}

\GG{180}{111111, 11112, 21111, 11121, 12111, 11211, 1122, 2211, 1212, 2121, 1221, 2112, 222, 1113, 3111, 1131, 1311, 123, 321, 132, 231, 213, 312}

\GG{181}{111111, 11112, 21111, 11121, 12111, 11211, 1122, 2211, 1212, 2121, 1221, 2112, 222, 1113, 3111, 1131, 1311, 123, 321, 132, 231, 213, 312, 33}

\GG{182}{111111, 11112, 21111, 11121, 12111, 11211, 1122, 2211, 1212, 2121, 1221, 2112, 222, 1113, 3111, 1131, 1311, 123, 321, 132, 231, 213, 312}

\GG{183}{111111, 11112, 21111, 11121, 12111, 11211, 1122, 2211, 1212, 2121, 1221, 2112, 222, 1113, 3111, 1131, 1311, 123, 321, 132, 231, 213, 312}

\GG{184}{111111, 11112, 21111, 11121, 12111, 11211, 1122, 2211, 1212, 2121, 1221, 2112, 222, 1113, 3111, 1131, 1311, 123, 321, 132, 231, 213, 312}

\GG{185}{111111, 11112, 21111, 11121, 12111, 11211, 1122, 2211, 1212, 2121, 1221, 2112, 222, 1113, 3111, 1131, 1311, 123, 321, 132, 231, 213, 312}

\GG{186}{111111, 11112, 21111, 11121, 12111, 11211, 1122, 2211, 1212, 2121, 1221, 2112, 222, 1113, 3111, 1131, 1311, 123, 321, 132, 231, 213, 312, 33}

\GG{187}{111111, 11112, 21111, 11121, 12111, 11211, 1122, 2211, 1212, 2121, 1221, 2112, 222, 1113, 3111, \underline{1131}, \underline{1311}, \underline{123}, \underline{321}, \underline{132}, \underline{231}, 213, 312}

\GG{188}{111111, 11112, 21111, 11121, 12111, 11211, 1122, 2211, 1212, 2121, 1221, 2112, 222, 1113, 3111, 1131, 1311, 123, 321, 132, 231, 213, 312, 33}

\GG{189}{111111, 11112, 21111, 11121, 12111, 11211, 1122, 2211, 1212, 2121, 1221, 2112, 222, 1113, 3111, 1131, 1311, 123, 321, 132, 231, 213, 312, 141}

\GG{190}{111111, 11112, 21111, 11121, 12111, 11211, 1122, 2211, 1212, 2121, 1221, 2112, 222, 1113, 3111, 1131, 1311, 123, 321, 132, 231, 213, 312, 33, 114, 411, 141, 24, 42}

\GG{191}{111111, 11112, 21111, 11121, 12111, 11211, 1122, 2211, 1212, 2121, 1221, 2112, 222, 1113, 3111, 1131, 1311, 123, 321, 132, 231, 213, 312, 114, 411, 141}

\GG{192}{111111, 11112, 21111, 11121, 12111, 11211, 1122, 2211, 1212, 2121, 1221, 2112, 222, 1113, 3111, 1131, 1311, 123, 321, 132, 231, 213, 312, 33}

\GG{193}{111111, 11112, 21111, 11121, 12111, 11211, 1122, 2211, 1212, 2121, 1221, 2112, 222, 1113, 3111, 1131, 1311, 123, 321, 132, 231, 213, 312}

\GG{194}{111111, 11112, 21111, 11121, 12111, 11211, 1122, 2211, 1212, 2121, 1221, 2112, 222, 1113, 3111, 1131, 1311, 123, 321, 132, 231, 213, 312, 33, 114, 411, 141, 24, 42}

\GG{195}{111111, 11112, 21111, 11121, 12111, 11211, 1122, 2211, 1212, 2121, 1221, 2112, 222, 1113, 3111, 1131, 1311, 123, 321, 132, 231, 213, 312, 33, 114, 411, 141, 24, 42}

\GG{196}{111111, 11112, 21111, 11121, 12111, 11211, 1122, 2211, 1212, 2121, 1221, 2112, 222, 1113, 3111, 1131, 1311, 123, 321, 132, 231, 213, 312, 33}

\GG{197}{111111, 11112, 21111, 11121, 12111, 11211, 1122, 2211, 1212, 2121, 1221, 2112, 222, 1113, 3111, 1131, 1311, 123, 321, 132, 231, 213, 312, 33, 141}

\GG{198}{111111, 11112, 21111, 11121, 12111, 11211, 1122, 2211, 1212, 2121, 1221, 2112, 222, 1113, 3111, 1131, 1311, 123, 321, 132, 231, 213, 312, 33}

\GG{199}{111111, 11112, 21111, 11121, 12111, 11211, 1122, 2211, 1212, 2121, 1221, 2112, 222, 1113, 3111, 1131, 1311, 123, 321, 132, 231, 213, 312, 33, 114, 411, 141, 24, 42}

\GG{200}{111111, 11112, 21111, 11121, 12111, 11211, 1122, 2211, 1212, 2121, 1221, 2112, 222, 1113, 3111, 1131, 1311, 123, 321, 132, 231, 213, 312, 33, 114, 411, 141, 24, 42}

\GG{201}{111111, 11112, 21111, 11121, 12111, 11211, 1122, 2211, 1212, 2121, 1221, 2112, 222, 1113, 3111, 1131, 1311, 123, 321, 132, 231, 213, 312, 33, 141}

\GG{202}{111111, 11112, 21111, 11121, 12111, 11211, 1122, 2211, 1212, 2121, 1221, 2112, 222, 1113, 3111, 1131, 1311, 123, 321, 132, 231, 213, 312, 33}

\GG{203}{111111, 11112, 21111, 11121, 12111, 11211, 1122, 2211, 1212, 2121, 1221, 2112, 222, 1113, 3111, 1131, 1311, 123, 321, 132, 231, 213, 312, 33, 114, 411, 141, 24, 42}

\GG{204}{111111, 11112, 21111, 11121, 12111, 11211, 1122, 2211, 1212, 2121, 1221, 2112, 222, 1113, 3111, 1131, 1311, 123, 321, 132, 231, 213, 312, 33, 114, 411, 141, 24, 42}

\GG{205}{111111, 11112, 21111, 11121, 12111, 11211, 1122, 2211, 1212, 2121, 1221, 2112, 222, 1113, 3111, 1131, 1311, 123, 321, 132, 231, 213, 312, 33, 114, 411, 141, 24, 42}

\GG{206}{111111, 11112, 21111, 11121, 12111, 11211, 1122, 2211, 1212, 2121, 1221, 2112, 222, 1113, 3111, 1131, 1311, 123, 321, 132, 231, 213, 312, 33, 114, 411, 141, 24, 42}

\GG{207}{111111, 11112, 21111, 11121, 12111, 11211, 1122, 2211, 1212, 2121, 1221, 2112, 222, 1113, 3111, 1131, 1311, 123, 321, 132, 231, 213, 312, 33, 114, 411, 141, 24, 42}

\GG{208}{111111, 11112, 21111, 11121, 12111, 11211, 1122, 2211, 1212, 2121, 1221, 2112, 222, 1113, 3111, 1131, 1311, 123, 321, 132, 231, 213, 312, 33, 114, 411, 141, 24, 42, 15, 51}

%% file: main.bbl
\begin{thebibliography}{99}

\bibitem{z}
AIM Minimum Rank -- Special Graphs Work Group, Zero forcing sets and the minimum rank of graphs, \emph{Linear Algebra and its Applications}, \textbf{428} (2008), 1628--1648.

\bibitem{z+}
Francesco Barioli, Wayne Barrett, Shaun Fallat, H. Tracy Hall, Leslie Hogben, Bryan Shader, P. van der Driessche, Hein van der Holst, Zero forcing parameters and minimum rank problems, \emph{Linear Algebra and its Applications} \textbf{433} (2010) 401–-411.

\bibitem{BF}
Francesco Barioli, Shaun Fallat, On the eigenvalues of generalized and double generalized stars, \emph{Linear and Multilinear Algebra} \textbf{53} (2005), 269--291.

\bibitem{atom}
Wayne Barrett, Steve Butler, H. Tracy Hall, John Sinkovic, Wasin So,  Colin Starr, Amy Yielding, Computing inertia sets using atoms, \emph{Linear Algebra and its Applications} \textbf{436} (2012), 4489--4502. 

\bibitem{BIRS}
Wayne Barrett, Steve Butler, Shaun Fallat, H. Tracy Hall, Jephian Lin, Leslie Hogben, Bryan Shader, Michael Young, The inverse eigenvalue problem of a graph: multiplicities and minors, \href{https://arxiv.org/abs/1708.00064}{arXiv:1708.00064}.

\bibitem{SSP}
Wayne Barrett, Shaun Fallat, H. Tracy Hall, Leslie Hogben, Jephian Lin, Bryan Shader, Generalizations of the strong Arnold property and the minimum number of distinct eigenvalues of a graph, Electronic Journal of Combinatorics, \textbf{24} (2017), P2.40, 28 pp.

\bibitem{inertia}
Wayne Barrett, H. Tracy Hall, Raphael Loewy, The inverse inertia problem for graphs:  Cut vertices, trees, and a counterexample, \emph{Linear Algebra and its Applications} \textbf{431} (2009), 1147--1191.

\bibitem{q}
Beth Bjorkman, Leslie Hogben, Scarlitte Ponce, Carolyn Reinhart, Theodore Tranel, Applications of analysis to the determination of the minimum number of distinct eigenvalues of a graph, \href{https://arxiv.org/abs/1708.01821}{arXiv:1708.01821}.

\bibitem{FH}
Shaun Fallat, Leslie Hogben, Minimum rank, maximum nullity, and zero forcing number of graphs, in \emph{Handbook of Linear Algebra, 2nd edition}, L. Hogben, ed., CRC Press, Boca Raton, 36 pp.

\bibitem{FF}
Ros\'ario Fernandes, C.M. da Fonseca, The inverse eigenvalue problem for Hermitian matrices whose graphs are cycles, \emph{Linear and Multilinear Algebra} \textbf{57} (2009), 673--682.

\bibitem{testSSP}
H. Tracy Hall, Code  in Sage for determining if  a  matrix satisfies Strong  Arnold  Property (SAP), Strong Spectral Property (SSP), or Strong Multiplicity Property (SMP). Sage worksheet available \url{https://sage.math.iastate.edu/home/pub/79/}.

\bibitem{H}
Leslie Hogben, Spectral graph theory and the inverse eigenvalue problem of a graph, \emph{Electronic Journal of Linear Algebra} \textbf{14} (2005), 12--31.

\bibitem{atlas}
Ronald C. Read, Robin J. Wilson, \emph{An Atlas of Graphs}, Oxford University Press, New York, 1998, xii+454 pp.

\end{thebibliography}
